 \documentclass{article} 
\usepackage[small,compact]{titlesec}
\usepackage{amsthm,amsmath,amsthm,amssymb,amsfonts,mathabx,mathrsfs,setspace,bbold,booktabs}
\usepackage{graphicx,caption,subfigure,epsfig,wrapfig,sidecap}
\usepackage[export]{adjustbox}

\usepackage{url,color,xcolor,verbatim,algorithmicx}
\usepackage[sort,compress]{cite}
\usepackage{diagbox}
\usepackage{algorithmicx}
\usepackage[noend]{algpseudocode}
\usepackage{bbm,bm}
\usepackage{arydshln}

\usepackage[ruled,boxed]{algorithm}

\newtheorem{theorem}{Theorem}

\newtheorem{corollary}[theorem]{Corollary}

\newtheorem{lemma}[theorem]{Lemma}
\newtheorem{proposition}[theorem]{Proposition}

\newtheorem{example}[theorem]{Example}

\def\R{\mathbb{R}}

\def\bY{\mathbf{Y}}

\def\mm#1{\boldsymbol{#1}}

\newcommand{\argmin}[1]{\underset{#1}{\operatorname{arg}\operatorname{min}}\;}

\newcommand{\argmax}[1]{\underset{#1}{\operatorname{arg}\operatorname{max}}\;}

\newcommand{\innerp}[1]{\langle{#1}\rangle}

\definecolor{mygrey}{gray}{0.75}

\newenvironment{rmat}{\left[\begin{array}{rrrrrrrrrrrrr}}{\end{array}\right]}
\newcommand\brm{\begin{rmat}}
\newcommand\erm{\end{rmat}}
\newenvironment{cmat}{\left[\begin{array}{ccccccccc}}{\end{array}\right]}
\newcommand\bcm{\begin{cmat}}
\newcommand\ecm{\end{cmat}}

\usepackage{soul}
\usepackage{array}

\begin{document}

\begin{center}

\textbf{\Large Robust First and Second-Order Differentiation for Regularized Optimal Transport}\\

\vspace{4mm}Xingjie Li\footnote{Department of Mathematics and Statistics, University of North Carolina, at Charlotte},  Fei Lu\footnote{ Department of Mathematics, Johns Hopkins University}, Molei Tao\footnote{School of Mathematics,
Georgia Institute of Technology}, Felix X.-F. Ye\footnote{Department of Mathematics \& Statistics,  University at Albany}\\
\end{center}

\begin{abstract}
Applications such as unbalanced and fully shuffled regression can be approached by optimizing regularized optimal transport (OT) distances, such as the entropic OT and Sinkhorn distances. A common approach for this optimization is to use a first-order optimizer, which requires the gradient of the OT distance. For faster convergence, one might also resort to a second-order optimizer, which additionally requires the Hessian. The computations of these derivatives are crucial for efficient and accurate optimization. However, they present significant challenges in terms of memory consumption and numerical instability, especially for large datasets and small regularization strengths.
We circumvent these issues by analytically computing the gradients for OT distances and the Hessian for the entropic OT distance, which was not previously used due to intricate tensor-wise calculations and the complex dependency on parameters within the bi-level loss function. 
Through analytical derivation and spectral analysis, we identify and resolve the numerical instability caused by the singularity and ill-posedness of a key linear system. Consequently, we achieve scalable and stable computation of the Hessian, enabling the implementation of the stochastic gradient descent (SGD)-Newton methods. 
Tests on shuffled regression examples demonstrate that the second stage of the SGD-Newton method converges orders of magnitude faster than the gradient descent-only method while achieving significantly more accurate parameter estimations.

\end{abstract}

\section{Introduction}

Optimal transport (OT) provides a powerful tool for finding a map between source and target distributions, especially when they are represented by ensemble samples without correspondence. Examples include shuffled regression \cite{pananjady2017linear,abid2017linear,hsu2017linear}, unlabeled sensing \cite{unnikrishnan2015unlabeled,unnikrishnan2018unlabeled,elhami2017unlabeled,zhang2020optimal}, homomorphic sensing\cite{tsakiris2018algebraic,tsakiris2019homomorphic}, regression with an unknown permutation \cite{li2021generalized}, or more broadly, as regression without correspondence \cite{rahimi07unsupervised,hsu2017linear,nejatbakhsh2019robust,xie2021hypergradient,AKLM22}.
 
The task is to find a parameterized function $ y=F( x; \theta)$ that maps ensemble of sources  $\mathbf{X} =\{\mathbf{x}_i\}_{i=1}^M \in \mathbb{R}^{M\times D}$ to targets $\mathbf{Y}^*=\{\mathbf{y}^*_j \}_{j=1}^N\in  \mathbb{R}^{N\times d}$ with probability weights $\mm \mu$ and $\mm \nu$. Here, $\mm \mu = (\mu_1,\ldots, \mu_M)^\top $ and  $\mm \nu = (\nu_1,\ldots, \nu_N)^\top $ satisfy $\mm \mu^\top  \mathbb{1}_M= \sum_{i=1}^M \mu_i=1, \mm\nu^\top \mathbb{1}_N = \sum_{j=1}^N \nu_j =1$ and $0<\mu_i,\nu_j<1$. Note that $M$ may not be equal to $N$, and the same applies to $D$ and $d$. The absence of a one-to-one correspondence between the source and target data samples makes classical supervised regression methods inapplicable.

The OT solution finds an optimal $\theta$ by minimizing a loss function $\mathcal{L}(\theta)$ between the image of the source data $\mathbf{Y}_\theta=F(\mathbf{X}; \theta)$ and the target data $\mathbf{Y}^*$, that is,
\begin{align} \label{framework_prob}
 \min_{\theta}  \mathcal{L}(\mathbf{C}, \mm\mu, \mm\nu),  
\end{align}
where the cost matrix $C_{ij} =c(\mathbf{y}_i(\theta), \mathbf{y}_j^*)$ and $c$ is a function of cost between $\mathbf{y}_i$ and $\mathbf{y}_j^*$. Throughout this study, we assume that $c(\mathbf{y},\mathbf{y}^*)$ is twice-differentiable, for instance, the squared Euclidean distance $c(\mathbf{y},\mathbf{y}^*)=\|\mathbf{y}-\mathbf{y}^*\|_2^2$. Several popular OT distances are candidates for the loss function, including the Wasserstein distance, the entropy-regularized OT (EOT) distance, and the Sinkhorn distance \cite{peyre2019computational}; see Section~\ref{subsec:loss_functions} for a brief review. Each of them leads to a bi-level optimization problem; for example, the EOT distance $\text{OT}_\epsilon(\mathbf{C},  \mm{\mu}, \mm\nu)$ leads to 
\begin{align*}
       &\min_\theta 
\, \, 
       \min_{\mm\Pi\in \mathbf{U}(\mm \mu, \mm \nu)}\sum_{i=1}^M\sum_{j=1}^N C_{ij} \Pi_{ij} + \epsilon \text{KL}( \mm\Pi, \mm{\mu} \otimes \mm\nu),  
\end{align*}
{\color{black}{where the convex polytope $ \mathbf{U}(\mm \mu, \mm \nu) = \{\mm\Pi\in\mathbb{R}^{M\times N}_{\ge 0}: \mm\Pi\mathbb{1}_N = \mm \mu, \mm\Pi^\top \mathbb{1}_M=\mm \nu\}$ is the set of bounded matrices with $M+N$ equality constraints.} }

A crucial component of the optimization process is computing the derivatives of the loss function with respect to $\theta$, and hence the derivatives of the OT distance with respect to data  
$\mathbf{Y}_\theta$. First-order optimization methods require the gradients of the OT distances. Danskin's theorem provides analytical gradients for the EOT distance \cite{bertsekas1997nonlinear, cuturi2014fast, genevay2018learning, feydy2019interpolating}, but it does not apply to the Sinkhorn distance. Hence, the generic implicit differentiation method \cite{luise2018differential} has been widely applied to OT distances \cite{cuturi2020supervised,xie2020differentiable, xie2021hypergradient,eisenberger2022unified, blondel2022efficient}. However, first-order optimization methods often converge slowly.

 To accelerate the convergence, a common strategy involves using stochastic gradient descent (SGD) initially, followed by Newton method's iterations, which necessitates computing the Hessian. Automatic differentiation and implicit differentiation are the two main methods for computing the Hessian \cite{blondel2022efficient,cuturi2022optimal,cuturi2020supervised}; see Section \ref{sec:prev} for a detailed discussion. However, both methods encounter significant challenges, such as memory shortages when the dataset is large and numerical instability due to singularity or ill-posedness, particularly when the entropy regularization strength $\epsilon$ is small. These issues impede the success of the SGD-Newton strategy.

We solve these issues by introducing analytical gradients for the OT distances and an analytical Hessian for the EOT. In particular, we achieve scalable and stable computation of the Hessian by using the analytical expression to locate and resolve the singularity or ill-posedness through spectral analysis. Our algorithm~\ref{Alg_Hess} significantly outperforms the automatic differentiation and implicit differentiation method in runtime and accuracy by orders of magnitudes; see Section \ref{sec:Hess_comp}. As a result, we enable the success of the SGD-Newton strategy for accelerating the bi-level optimization, as we demonstrate on parameter estimation for shuffled regression of mixed Gaussian and 3D Point Clouds Registration of MobilNet10 dataset \cite{qi2017pointnet,katageri2024metric}; see Section \ref{sec:Application}.
 
 The key in our derivation is the linear system for the optimal dual potentials, e.g., \eqref{linear-sys-k} or \eqref{linear-sys}, which is inspired by the implicit differentiation in \cite{blondel2022efficient} and the second-order Fr\'{e}chet derivative of the Sinkhorn divergence loss under the Wasserstein metric in \cite{shen2020sinkhorn}.
 Emerging from the implicit differentiation, this linear system facilitates efficient computation of the gradient of the OT distances as well as the Hessian of the EOT distance. In particular, when used together with Eq.\eqref{key-observation} from the marginal constraints, it bridges implicit differentiation and Danskin's theorem in the context of EOT distance. 
 
Furthermore, we provide a comprehensive spectral analysis for the linear system for the dual potentials through the matrix 
 \begin{equation}\label{eq:H}
\mathbf{H}:=  	    \bcm
        \text{diag}( \mm\Pi \mathbb{1}_N) & \mm\Pi \\ (\mm\Pi)^\top & \text{diag}(\mm\Pi^\top  \mathbb{1}_M)
    \ecm \in \R^{(M+N)\times (M+N)}, 
 \end{equation}
 where $\mm\Pi$ is the coupling matrix. {\color{black}{When $\mm\Pi^*$ is the optimal coupling matrix, we denote this matrix \eqref{eq:H} as $\mathbf{H}^*$.}}
 We show that when $\mm\Pi$ has positive entries, $\mathbf{H}$ has zero as a simple eigenvalue, and its effective condition number (i.e., the ratio of the largest and smallest positive eigenvalues) has upper and lower bounds depending on the spectral gap of $\mm\Pi^\top \mm\Pi$. In particular, we construct an example showing that $\mathbf{H}$ can be severely ill-conditioned with the smallest positive eigenvalue at the order of $O(e^{-\frac{1}{\epsilon}})$ when $\epsilon$ is small, or $O(\frac{1}{N})$ when $N$ is large. Thus, when solving a linear system with $\mathbf{H}$, proper regularization is crucial.

Our main contributions are threefold. 

\begin{itemize}
\item \emph{Analytical derivatives and spectral analysis.} We derive analytical gradients with respect to the data $\mathbf{Y}$ for EOT and Sinkhorn distances and Hessian for the EOT distance in Section~\ref{sec:derivative}--\ref{subsec:Hessian_EOT}. The spectral analysis in Section \ref{sec:spectrum} helps us understand and resolve the numerical instability issue via a proper regularization in Section \ref{subsec:TrunSVD}.

\item \emph{Fast stable computation of Hessian.} Our algorithm enables a stable, memory-efficient, and fast computation of the Hessian, significantly outperforming other state-of-the-art methods in runtime and accuracy by orders of magnitudes; see Section \ref{sec:Hess_comp}.

\item \emph{Enabling SGD-Newton for shuffled regression.} With the robust computation of the Hessian, we are able to apply the SGD-Newton method to shuffled regression problems in Section \ref{sec:Application}, significantly accelerating the optimization process.
\end{itemize}

\subsection{Outline}
This work is organized as follows.  Section~\ref{sec:background} reviews the various (OT) distances and the Sinkhorn algorithm. Section \ref{sec:derivatives_EOT} is devoted to the analytical and numerical computation of the gradients and Hessian, leading to an algorithm with proper regularization. Then, we analyze the spectrum of the matrix $\mathbf{H}$ in Section \ref{sec:spectrum}. 
 In Section~\ref{sec:Hess_comp}, we examine the efficiency and accuracy of Hessian computation using a benchmark example and compare the results with other approaches.  Then we apply the proposed SGD-Newton approach to applications in Section~\ref{sec:Application}, including the parameter estimation for shuffled regression of mixed Gaussian and 3D Point Clouds Registration of MobilNet10 dataset.

\section{Optimal Transport Loss and Sinkhorn Algorithm}\label{sec:background} 

Ideally, we could find 
$\theta^*$ by minimizing the optimal transport distance between the parameterized source data $\mathbf{Y}_\theta$ and the target data $\mathbf{Y}^*$. 
We will first review some classical results in computational optimal transport \cite{peyre2019computational} in this section. In the section, we ignore $\theta$ in the notation unless noted otherwise. 
\subsection{ Optimal Transport Loss Functions}\label{subsec:loss_functions}
\subsubsection{Wasserstein-2 Metric}
 One popular choice is to use the Wasserstein-2 metric as the optimal transport loss, equivalently, $ \mathcal{L}( \mathbf{C}, \mm \mu,  \mm \nu)=W_2^2(\mathbf{C}, \mm \mu,\mm \nu)$. 
To calculate the Wasserstein-2 metric, one has to solve a {\color{black}{constrained optimization}} problem,   
\begin{align}\label{Wasserstein-2}
    &W_2^2(\mathbf{C}, {\mm \mu},  \mm \nu) :=  \min_{\mm\Pi\in \mathbf{U}(\mm \mu, \mm \nu)}\sum_{i=1}^M\sum_{j=1}^N C_{ij} \Pi_{ij}, \quad  C_{ij}=  c(\mathbf{y}_i,\mathbf{y}_j^*) , 
\end{align}
where  $ c(\mathbf{y}_i,\mathbf{y}_j^*)  =\| \mathbf{y}_i-\mathbf{y}_j^*\|_2^2 $ is the cost of transport, and the \emph{coupling matrix} $\mm\Pi \in [0,1]^{M\times N}$ is the transport plan from the parameterized source data $\mathbf{Y}$ to the target data $\mathbf{Y}^*$.  To solve the constrained optimization 
via linear programming, the computational complexity is $O((N + M)NM \log(N+ M))$\cite{peyre2019computational}, which is very expensive when $N,M$ are large.  To overcome this issue, one often regularizes the objective function. Common regularization includes the EOT 
distance, the Sinkhorn distance, which we briefly review below.

\subsubsection{Entropy-regularized OT (EOT) Distance}

The EOT distance is the Wasserstein-2 loss plus the relative entropy
 between two measures:
\begin{align}\label{reg_OT_problem}
       &\text{OT}_\epsilon(\mathbf{C},  \mm{\mu}, \mm\nu):= \min_{\mm\Pi\in \mathbf{U}(\mm \mu, \mm \nu)}\sum_{i=1}^M\sum_{j=1}^N C_{ij} \Pi_{ij} + \epsilon \text{KL}(\mm\Pi, \mm{\mu} \otimes \mm\nu), 
\end{align}
where the {\it relative entropy} between the coupling matrix $\mm\Pi$ and the outer product $\mm\mu\bigotimes \mm\nu$ is 
$\text{KL}(\mm\Pi, \mm{\mu} \otimes \mm\nu): = \sum_{i=1}^M \sum_{j=1}^N \big(\Pi_{ij} \log \frac{\Pi_{ij}}{\mu_i \nu_j}{\color{black}-\Pi_{ij}+\mu_i\nu_j}\big)$ \cite{feydy2020geometric}.  
This regularization drastically simplifies the study of the dual problem and further leads to the Sinkhorn algorithm for its unique numerical solution \cite{peyre2019computational}.
As $\epsilon$ goes to 0, EOT converges to the 
Wasserstein-2 distance
at the rate of $\epsilon$ \cite{luise2018differential}.

\subsubsection{Sinkhorn Distance}
Another candidate for the regularized OT loss is called the Sinkhorn distance,  $\widetilde{\text{OT}}_\epsilon(\mathbf{C}, \mm{\mu}, \mm\nu)$, 
\begin{align} \label{sinkhorn_dist}
     & \widetilde{\text{OT}}_\epsilon(\mathbf{C}, \mm{\mu}, \mm\nu) := \sum_{ij} C_{ij} \Pi^*_{ij},   \text{ with }  \mm\Pi^*=\argmin{\mm\Pi\in \mathbf{U}(\mm \mu, \mm \nu)}\sum_{ij}C_{ij}\Pi_{ij} + \epsilon \text{KL}(\mm\Pi, \mm{\mu}\otimes \mm\nu).
\end{align}
Sinkhorn distance eliminates the contribution of the entropy regularization term from $\text{OT}_\epsilon(\mathbf{C}, \mm{\mu}, \mm\nu)$ to the total loss $\mathcal{L}$ after the transport plan $\mm\Pi^*$ has been obtained.  It gives even better approximation results and converges to the Wasserstein-2 distance exponentially fast, that is, we have 
$\left|\widetilde{\text{OT}}_\epsilon(\mathbf{C}, \mm{\mu}, \mm\nu) - W_2^2(\mathbf{C}, \mm{\mu}, \mm\nu)\right|\le c\exp(-1/\epsilon)$ \cite{luise2018differential}.

\subsection{ Sinkhorn Algorithm}\label{sec:sinkhorn}
From now on, we firstly choose the EOT cost as the loss function, i.e., $\mathcal{L}(\mathbf{C}, \mm{\mu}, \mm\nu)=\text{OT}_\epsilon (\mathbf{C}, \mm{\mu}, \mm\nu)$ though we will discuss the others later. 
The computation of this quantity, i.e., the constrained optimization \eqref{reg_OT_problem} is solved by the well-known Sinkhorn algorithm. As a preliminary, we briefly review this algorithm.  
One can introduce two slack variables, known as the dual potentials, $\mathbf{f}\in \mathbb{R}^M, \mathbf{g}\in \mathbb{R}^N$, for each marginal constraint of \eqref{reg_OT_problem}.  
{\color{black}Also because of the relative entropy term $\log\frac{\Pi_{ij}}{\mu_i\nu_j}$, $\Pi_{ij}$ will never be negative and we can drop non-negativity constraints $\Pi_{ij}\ge 0$.} As a result, the corresponding 
{\color{black}Lagrangian} is 
\begin{equation*}
\begin{split}
    L(\mathbf{C}, \mm\Pi, \mathbf{f},  \mathbf{g})&= \sum_{ij}\left(C_{ij}\Pi_{ij} + \epsilon\Pi_{ij}\big(\log \frac{\Pi_{ij}}{\mu_i\nu_j}{\color{black}-1}\big)- f_i (\Pi_{ij}-\mu_i) - g_j( \Pi_{ij}-\nu_j)\right) \color{black}{+\epsilon}\, . 
    \end{split}
\end{equation*}

\noindent
 {\color{black} The strong duality \cite{bertsimas1997introduction,clason2021entropic} holds, $\text{OT}_\epsilon(\mathbf{C}, \mm \mu, \mm \nu)=\max_{\mathbf{f}, \mathbf{g}}\min_{\mm \Pi}L(\mathbf{C}, \mm \Pi, \mathbf{f}, \mathbf{g})$. Also, the first-order optimality conditions of the inner problem
 yield the expression of the optimal coupling matrix $\mm\Pi^*$ as the function of $\mathbf{f}$ and $\mathbf{g}$, $\Pi_{ij}^* = \mu_i\nu_j \exp\left(\frac{-C_{ij}+f_i+g_j}{\epsilon}\right)$. It ensures that the optimal coupling matrix $\mm\Pi^*$ is entry-wise positive.  The first-order optimality conditions of the outer problem yield the marginal constraints, 
\begin{align}\label{constraint}
    \sum_j \mu_i\nu_j \exp\left(\frac{-C_{ij}+f^*_i+g^*_j}{\epsilon}\right) = \mu_i, \ \ \sum_i \mu_i\nu_j \exp\left(\frac{-C_{ij}+f^*_i+g^*_j}{\epsilon}\right) = \nu_j.
\end{align}
}

An intuitive scheme to solve these nonlinear equations \eqref{constraint} for optimal dual variable $\mathbf{f}^*$ and $\mathbf{g}^*$ is to alternatively rescale rows and columns of the Gibbs kernel to satisfy the marginal
constraint, which is called \emph{Sinkhorn algorithm}\cite{cuturi2013sinkhorn,sinkhorn1967concerning, sinkhorn1964relationship,marshall1968scaling}. Numerically, however, this computation becomes unstable when $\epsilon$ is small. 
The stable Sinkhorn iteration is thus performed in the log-domain \cite{cuturi2013sinkhorn, peyre2019computational}, 
\begin{align*}
    &\mathbf{f}^{(l+1)} = \epsilon \log \mm \mu -\epsilon \log \left(\mathbf{K} \exp(\mathbf{g}^{(l)}/\epsilon)\right),\; \mathbf{g}^{(l+1)} = \epsilon \log \mm \nu - \epsilon \log \left(\mathbf{K}^\top \exp(\mathbf{f}^{(l+1)}/\epsilon)\right)
\end{align*}
with the initial vector to be $\mathbf{g}^{(0)}=\mathbb{0}_N$ and the Gibbs kernel $\mathbf{K} = \exp\left(-\frac{\mathbf{C}}{\epsilon}\right)$. As $l$ goes to $+\infty$, both converge to $\mathbf{f}^*$ and $\mathbf{g}^*$. In practice, the iteration stops when the 1-norm of marginal violation is within the threshold value. 

\noindent
{\color{black}{Finally we get the optimal coupling matrix $\mm\Pi^*$ and  EOT distance $ \text{OT}_\epsilon(\mathbf{C},  \mm\mu, \mm\nu)$,
\begin{align}\label{KOT_dual}     &\Pi_{ij}^* = \mu_i\nu_j \exp\left(\frac{-C_{ij}+f^*_i+g^*_j}{\epsilon}\right), \ \  \text{OT}_\epsilon(\mathbf{C},  \mm{\mu}, \mm\nu)
     = \bcm \mm \mu^\top & \mm \nu^\top \ecm \bcm  \mathbf{f}^* \\ 
     \mathbf{g}^* \ecm.
\end{align}}} 
Overall the computational complexity of EOT to achieve $\tau$-approximate of the unregularized OT problem is $O(N^2\log(N)\tau^{-3})$ when $M=N$ \cite{altschuler2017near, peyre2019computational}, {\color{black}{which is a significant improvement}} to the linear programming of Wasserstein-2 metric. 

Similarly,
$\widetilde{\text{OT}}_\epsilon(\mathbf{C},  \mm{\mu}, \mm\nu)$  
gives
\begin{align}\label{sink_distance}
    &\widetilde{\text{OT}}_\epsilon(\mathbf{C},  \mm{\mu}, \mm\nu) = \sum_{ij}C_{ij}\Pi^*_{ij}= \sum_{ij}\mu_i\nu_j C_{ij}e^{\frac{-C_{ij} +f_i^* +g_j^*}{\epsilon}} \,. 
\end{align}

\section{Differentiation of Loss Functions}\label{sec:derivatives_EOT}
In this section, we introduce robust computations for the gradients of these regularized OT loss functions and for the Hessian of the EOT distance. 
Our computations show that costly backward propagation can be avoided even for the Hessian. 

We recall the problem setup \eqref{framework_prob} and first study the analytic form of the gradient and hessian of EOT distance with respect to the parameter $\theta$, which read
\begin{align}\label{deriv_theta}
&\frac{\partial \text{OT}_\epsilon(\mathbf{C}_\theta,  \mm\mu, \mm\nu)}{\partial \theta_i} = \sum_{k}\frac{\partial {\mathbf{y}}_k}{\partial \theta_i} \frac{\partial \text{OT}_\epsilon(\mm 
 C_\theta,  \mm\mu, \mm \nu)}{\partial  {\mathbf{y}}_k},  \\ 
   & \frac{\partial ^2\mathrm{OT}_{\epsilon}(\mathbf{C}_\theta, \mm{\mu}, \mm\nu)}{\partial \theta_i \partial \theta_j} = \sum_{s,k} \frac{\partial {\mathbf{y}}_s}{\partial \theta_i} \frac{\partial ^2\mathrm{OT}_{\epsilon}(\mathbf{C}_\theta, \mm{\mu}, \mm\nu)}{\partial {\mathbf{y}}_s \partial {\mathbf{y}}_k} \left(\frac{\partial {\mathbf{y}}_k}{\partial \theta_j}\right)^\top\nonumber  + \sum_{k} \frac{\partial ^2{\mathbf{y}}_k}{\partial \theta_i \partial \theta_j}\frac{\partial \mathrm{OT}_\epsilon(\mathbf{C}_\theta, \mm{\mu}, \mm\nu)}{\partial {\mathbf{y}}_k}. 
  \end{align}  
The key step is to find the explicit expression for the first and second derivatives  with respect to the source data 
$\frac{\partial \text{OT}_\epsilon(\mathbf{C}_\theta,  \mm \mu,  \mm\nu)}{\partial  {\mathbf{y}}_k}\in \mathbb{R}^d$, 
$\frac{\partial ^2\mathrm{OT}_{\epsilon}(\mathbf{C}_\theta, \mm\mu,  \mm\nu)}{\partial {\mathbf{y}}_s \partial {\mathbf{y}}_k}\in \mathbb{R}^{d\times d} $.
{\color{black}{To provide a full picture, the gradients of EOT distance and Sinkhorn distance with respective to the measures $\mm \mu$ or  $\mm \nu$ have been studied in \cite{luise2018differential,peyre2019computational}.} }

\subsection{Previous methods computing  first- and second-order derivatives}\label{sec:prev}
There are two main approaches to the analytical expression of the first-order derivatives of regularized OT distances: the Danskin approach and the implicit differentiation approach. 
The Danskin approach computes the gradient of EOT, based on applying the Danskin's theorem to the dual function \cite{bertsekas1997nonlinear, cuturi2014fast, genevay2018learning, feydy2019interpolating}.   
{\color{black}{
\begin{theorem}[Danskin's theorem]
    Suppose $\phi(\mathbf{x}, \mathbf{z})$ is a continuous function of two arguments, $\phi: \mathbb{R}^n\times Z\rightarrow \mathbb{R}$, where $Z\subset \mathbb{R}^m$ is a compact set and 
    define $f(\mathbf{x})=\max_{\mathbf{z}\in Z}\phi(\mathbf{x}, \mathbf{z})$ and the set of maximizing points $\mathbf{z}^*(\mathbf{x})=\argmax{\mathbf{z}\in Z}\phi(\mathbf{x}, \mathbf{z})$. If $\mathbf{z}^*(\mathbf{x})$ consists of a single element, the gradient of $f(\mathbf{x})$ is given by $\frac{\partial f}{\partial \mathbf{x}}=\frac{\partial \phi(\mathbf{x}, \mathbf{z}^*(\mathbf{x}))}{\partial \mathbf{x}}$.
\end{theorem}

 Recall the outer problem $\text{OT}_\epsilon(\mathbf{C}, \mm \mu, \mm \nu)= \max_{ \mathbf{f}, \mathbf{g}}L(\mathbf{C},\mm\Pi^*(\mathbf{f}, \mathbf{g}), \mathbf{f}, \mathbf{g}) $, 
 \begin{align*}
     \text{OT}_\epsilon(\mathbf{C}, \mm \mu, \mm \nu)= \max_{ \mathbf{f}, \mathbf{g}} \underbrace{\sum_i f_i \mu_i +\sum_j g_j\nu_j - \sum_{ij}\epsilon \mu_i\nu_j \exp\left(\frac{-C_{ij}+f_i+g_j}{\epsilon}\right) +\epsilon}_{\phi(\mathbf{Y},\mathbf{f},\mathbf{g})}.
 \end{align*}
If $(\mathbf{f}^*, \mathbf{g}^*)$ is a single element, applying Danskin's theorem to  $\phi(\mathbf{Y},\mathbf{f},\mathbf{g})$, we obtain 
\begin{align}
\label{eq:danskin}
    \frac{\partial \text{OT}_\epsilon(\mathbf{C}, \mm \mu, \mm \nu)}{\partial \mathbf{y}_k} &=\sum_{ij}  \mu_i\nu_j \exp\left(\frac{-C_{ij}+f^*_i+g^*_j}{\epsilon}\right) \frac{\partial C_{ij}}{\partial \mathbf{y}_k}=\sum_{j}\frac{\partial C_{kj}}{\partial \mathbf{y}_k}\Pi^*_{kj}.
\end{align}
The last equality is due to $\frac{\partial C_{ij}}{\partial \mathbf{y}_k}=\frac{\partial C_{kj}}{\partial \mathbf{y}_k}\delta_{ki}$.  }}
Once $\mm\Pi^*$ is obtained from the Sinkhorn algorithm, the first-order derivative with respect to the source dataset $\mathbf{Y}$ immediately follows without additional computational cost. {\color{black}{Unfortunately Danskin's approach  
  doesn't work on computing the gradient of the Sinkhorn distances} $\widetilde{\text{OT}}_\epsilon(\mathbf{C}_\theta, \mm{\mu}, \mm\nu)$ because it is not of the form  $\max_{\mathbf{f}, \mathbf{g}} \phi(\mathbf{C}, \mathbf{f}, \mathbf{g})$}. 
An implicit differentiation approach is thus introduced in \cite{luise2018differential}. It implicitly differentiates the associated marginal constraints to derive a large linear system, which is solved by the conjugate gradient solver in \textit{lineax} \cite{rader2023lineax}. It applies to different regularized OT problems \cite{cuturi2020supervised,xie2020differentiable, xie2021hypergradient, eisenberger2022unified, blondel2022efficient}.

For computing the Hessian of the OT distances, up to our knowledge, there is no direct analytical expression before the current work. {\color{black}{Similarly, Danskin's approach doesn't work on computing the Hessian of EOT distance. }}
There are two approaches suggested by \textit{OTT}\cite{cuturi2022optimal}. The first approach \textit{unrolls} the Sinkhorn iterations and use the \textit{JAX} in-build tools to handle autodiff via backward propagation and computational graph. The second approach \textit{implicitly} differentiates the optimal solutions computed by \textit{OTT}. 
The implicit differentiation approach involves differentiating the solution of an ill-conditioned linear system with the custom differentiation rules \cite{cuturi2022optimal}, hence, regularization techniques, such as preconditioning the marginal constraints \cite{cuturi2020supervised, ye2024enhancing} and the ridge regularization, have been introduced to try to resolve this issue. 
However, as we will demonstrate in Section~\ref{sec:Hess_comp}, both approaches still encounter two major challenges: (i) memory shortages when the dataset 
is large, and (ii) numerical instability due to singularity and the ill-posed nature of the linear system, particularly when 
$\epsilon$ is small.

The key in our derivation is the 
linear system for the optimal dual potentials, e.g., \eqref{linear-sys-k} or \eqref{linear-sys}, which {\color{black}{arises from applying implicit differentiation}}. It facilitates efficient computation of the gradient of the Sinkhorn distance as well as the Hessian of the EOT distance. Additionally, together with \eqref{key-observation}, it bridges implicit differentiation and Danskin's theorem in the context of the EOT distance.


\subsection{Analytical computation of the gradients}\label{sec:derivative}
 
We first review the implicit differentiation approach to the gradient of the EOT distance $\text{OT}_{\epsilon}(\mathbf{C}, \mm \mu,  \mm \nu)$ with respect to source data $\mathbf{Y}$ and re-derive the result of Danskin's approach in \eqref{eq:danskin} through the key observation.  
{\color{black}{With the same technique}}, we further provide a novel numerical method to efficiently compute the derivative of Sinkhorn distance $\widetilde{\mathrm{OT}}_\epsilon (\mathbf{C},  \mm\mu,  \mm\nu)$.

\bigskip
\noindent
\textbf{Gradient of $\text{OT}_\epsilon(\mathbf{C},  \mm \mu,  \mm \nu)$.} 
We first consider the gradient of the EOT distance with respect to the source data $\mathbf{y}_k$
\begin{subequations}
\label{derivatives}
  \begin{equation} \label{deriv_OneSample}
  \frac{\partial \text{OT}_\epsilon(\mathbf{C}, \mm \mu,  \mm\nu)}{\partial  \mathbf{y}_k}= \sum_{i=1}^M \mu_i \frac{\partial f_i^*}{\partial \mathbf{y}_k} + \sum_{j=1}^N  \nu_j\frac{\partial g^*_j}{\partial \mathbf{y}_k} = \bcm \mm\mu^\top & \mm\nu^\top \ecm
\bcm \frac{\partial \mathbf{f}^*}{\partial \mathbf{y}_k} \\ \frac{\partial \mathbf{g}^*}{\partial \mathbf{y}_k} \ecm,
    \end{equation}
    where $ \frac{\partial \mathbf{f}^*}{\partial \mathbf{y}_k} = \left( \frac{\partial f^*_1}{\partial \mathbf{y}_k}, \ldots, \frac{\partial f^*_M}{\partial \mathbf{y}_k} \right)^\top \in \R^{M\times d}$, and similarly $ \frac{\partial \mathbf{g}^*}{\partial \mathbf{y}_k} \in \R^{N\times d}$. To simplify the notation, we denote
$\left(\frac{d\mathbf{f}^*}{d\mathbf{Y}}\right)_{iks}=\frac{\partial f^*_i}{\partial (\mathbf{y}_k)_s}$ and $\left(\frac{d\mathbf{g}^*}{d\mathbf{Y}}\right)_{jks}=\frac{\partial g^*_j}{\partial (\mathbf{y}_k)_s}$, so $\frac{d\mathbf{f}^*}{d\mathbf{Y}}\in \mathbb{R}^{M\times M \times d}$ and $\frac{d\mathbf{g}^*}{d\mathbf{Y}}\in \mathbb{R}^{N\times M \times d}$.
    
We write the gradient in the vector form
    \begin{equation} \label{deriv_AllSample}
    \frac{d\text{OT}_\epsilon(\mathbf{C},  \mm \mu, \mm \nu)}{d \mathbf{Y}}
    =  \mm \mu^\top \frac{d\mathbf{f}^*}{d\mathbf{Y}} + \mm \nu^\top \frac{d \mathbf{g}^*}{d \mathbf{Y}}   = \bcm \mm \mu^\top & \mm \nu^\top \ecm   \begin{cmat}
    \frac{d \mathbf{f}^*}{d \mathbf{Y} } \\  \frac{d \mathbf{g}^*}{d \mathbf{Y} } 
    \end{cmat} \in \R^{ M\times d}.
\end{equation}
\end{subequations}

\begin{theorem}\label{thm:derivative} 
 The derivative of EOT distance in \eqref{KOT_dual} with respect to source data $\mathbf{y}_k$, as in \eqref{deriv_OneSample}, is given by    
    \begin{align}\label{derivate_data}
    \frac{\partial \mathrm{OT}_\epsilon (\mathbf{C}, \mm{\mu},  \mm\nu)}{\partial  \mathbf{y}_k}= \sum_{j=1}^N \frac{\partial C_{kj}}{\partial {\mathbf{y}_k}}\Pi^*_{kj},\quad k=1,\dots, M.
    \end{align}
   In vector form, the derivative of  $\mathrm{OT}_\epsilon (\mathbf{C},  \mm \mu, \mm \nu)$ with respect to whole source data $\mathbf{Y}$ is  
 \begin{align} \label{derivate_data_Y}
     \frac{d\mathrm{OT}_\epsilon (\mathbf{C}, \mm\mu,  \mm\nu)}{d\mathbf{Y}}=  \mathcal{B}\cdot \mathbb{1}_N, 
 \end{align}
where $\mathcal{B} \in \mathbb{R}^{M\times d\times N }$ is a tensor with entries $ \mathcal{B}_{ks j} = \frac{\partial C_{kj}}{\partial (\mathbf{y}_{k})_s}\Pi^*_{kj}$, 
and $\mathcal{B}\cdot\mathbb{1}_N \in \mathbb{R}^{M\times d}$ is the dot product which is the summation of the third index such that the $k$-th column is $(\mathcal{B}\cdot\mathbb{1}_N)_k = \sum_j\frac{\partial C_{kj}}{\partial \mathbf{y}_k}\Pi^*_{kj}$ for $1\leq k\leq M$. 
\end{theorem}
\begin{proof} 
The main task is to find $ \frac{\partial \mathbf{f}^*}{\partial \mathbf{y}_k}$ and $\frac{\partial \mathbf{g}^*}{\partial \mathbf{y}_k}$ in \eqref{deriv_OneSample} and the marginal probability conditions. 
    The partial derivative of entries of optimal coupling matrix $\mm\Pi^*$ with respect to $\mathbf{y}_k$ is 
    \begin{align}\label{EOT_Pi_derivative}
    \frac{\partial \Pi^*_{ij}}{\partial \mathbf{y}_k} &= \frac{\Pi^*_{ij} }{\epsilon} \left(-\frac{\partial C_{ij}}{\partial \mathbf{y}_k}  + \frac{\partial f_i^*}{\partial \mathbf{y}_k} + \frac{\partial g_j^*}{\partial \mathbf{y}_k}\right) = \frac{\Pi^*_{ij} }{\epsilon} \left(-\frac{\partial C_{kj}}{\partial \mathbf{y}_k}\delta_{ik}  + \frac{\partial f_i^*}{\partial \mathbf{y}_k} + \frac{\partial g_j^*}{\partial \mathbf{y}_k}\right),
    \end{align}
    where $\delta_{ik}$ is a Kronecker delta function. 
    We observe that with the marginal 
probability conditions $\sum_{j=1}^N \Pi^*_{ij}=\mu_i$ and $\sum_{i=1}^M \Pi^*_{ij}=\nu_j$, and by taking the partial derivative $\frac{\partial}{\partial \mathbf{y}_k}$ on both sides of these marginal constraints, we can get
    \begin{align*}
        & 0= \sum_{j=1}^N \frac{\partial \Pi^*_{ij}}{\partial \mathbf{y}_k} 
        =  \frac{  \mu_i}{\epsilon} \frac{\partial f_i^*}{\partial \mathbf{y}_k} -  \frac{1}{\epsilon}  \left[ \sum_{j=1}^N \left(\frac{\partial C_{kj}}{\partial \mathbf{y}_k}\delta_{ik} - \frac{\partial g^*_j}{\partial \mathbf{y}_k}\right) \Pi^*_{ij}\right], \\ 
     &  0= \sum_{i=1}^M \frac{\partial \Pi^*_{ij}}{\partial {\mathbf{y}_k}} 
     =  \frac{ \nu_j }{\epsilon} \frac{\partial g_j^*}{\partial \mathbf{y}_k} -  \frac{1}{\epsilon}  \left[ \sum_{i=1}^M \left(\frac{\partial C_{kj}}{\partial \mathbf{y}_k}\delta_{ik} -\frac{\partial f^*_i}{\partial \mathbf{y}_k}\right) \Pi^*_{ij}\right]. 
    \end{align*}
    Hence we have
 \begin{align*}
     &\mu_i \frac{\partial f_i^*}{\partial \mathbf{y}_k} +\sum_{j=1}^N \frac{\partial g^*_j}{\partial \mathbf{y}_k} \Pi^*_{ij} = \sum_{j=1}^N \frac{\partial C_{kj}}{\partial \mathbf{y}_k} \Pi^*_{ij}\delta_{ik} \in \R^d, \quad \forall 1\leq i\leq M, \\
     &  \sum_{i=1}^M \frac{\partial f^*_i}{\partial \mathbf{y}_k} \Pi^*_{ij} + \nu_j \frac{\partial g_j^*}{\partial \mathbf{y}_k} =   \frac{\partial C_{kj}}{\partial \mathbf{y}_k} \Pi^*_{kj}\in \R^d, \quad \forall 1\leq j\leq N.
\end{align*} 
The above linear system can be written in matrix form, 
\begin{align}  \label{linear-sys-k}
\bcm
        \text{diag}( \mm \mu) & \mm\Pi^* \\ (\mm\Pi^*)^\top & \text{diag}(\mm\nu)
    \ecm \bcm
    \frac{\partial \mm{f}^*}{\partial \mathbf{y}_k } \\  \frac{\partial \mm{g}^*}{\partial  \mathbf{y}_k } 
    \ecm 
    =  \bcm 
    \mathbf{e}_k  (\mathbf{B}_k\mathbb{1}_N)^\top \\
    \mathbf{B}_k^\top\ecm,
\end{align}
where the left-hand-side matrix is 
$\mathbf{H}^*$ defined in \eqref{eq:H}, 
 $\mathbf{e}_k$ is the $k$-th standard column basis of $\R^M$ and  the matrix $\mathbf{B}_k \in \mathbb{R}^{d\times N}$ is 
$(\mathbf{B}_k)_{sj}=\frac{\partial C_{kj}}{\partial (\mathbf{y}_k)_s} \Pi^*_{kj}$.
{\color{black}{Notice that $ \bcm \mathbb{1}_{M}^\top & -\mathbb{1}_N^\top\ecm \bcm 
    \mathbf{e}_k  (\mathbf{B}_k\mathbb{1}_N)^\top \\ \mathbf{B}_k^\top\ecm=\mathbb{0}_{d} $, and  $\text{span}\left(\bcm \mathbb{1}_{M}\\-\mathbb{1}_N\ecm\right)$ is  the kernel of $\mathbf{H}^*$ as proved in Theorem~\ref{prop:H_general_Pi},
so $\bcm 
    \mathbf{e}_k  (\mathbf{B}_k\mathbb{1}_N)^\top \\ \mathbf{B}_k^\top\ecm$ is in the span of $\mathbf{H}^*$,
therefore the linear system above always has a solution. }}

Instead of solving the above linear system to evaluate \eqref{deriv_OneSample}, we observe from the marginal constraints that 
\begin{align} \label{key-observation}
  \bcm \frac{1}{M}\mathbb{1}^\top_M & \frac{1}{N} \mathbb{1}^\top_N \ecm  
     \mathbf{H}^*
= \left(\frac{1}{M}+\frac{1}{N}\right) \bcm \mm \mu^\top & \mm\nu^\top\ecm. 
\end{align}
Then, multiplying both sides of \eqref{linear-sys-k} by $\bcm \frac{1}{M}\mathbb{1}^\top_M & \frac{1}{N} \mathbb{1}^\top_N \ecm$, we obtain 
\begin{align*}
    \bcm \frac{1}{M}\mathbb{1}^\top_M & \frac{1}{N} \mathbb{1}^\top_N \ecm  
    \mathbf{H}^*
     \bcm
    \frac{d \mathbf{f}^*}{d  \mathbf{y}_k } \\  \frac{d \mathbf{g}^*}{d  \mathbf{y}_k }     \ecm &
    = 
    \bcm \frac{1}{M}\mathbb{1}^\top_M & \frac{1}{N} \mathbb{1}^\top_N \ecm \bcm 
    \mathbf{e}_k  (\mathbf{B}_k\mathbb{1}_N)^\top \\
    \mathbf{B}_k^\top\ecm.
\end{align*}
Hence, the derivative of $\mathrm{OT}_\epsilon (\mathbf{C},  \mm\mu,  \mm\nu)$ with respect to source data $\mathbf{y}_k$ is 
\begin{align*}
       &\frac{d\mathrm{OT}_\epsilon (\mathbf{C}, \mm\mu, \mm \nu)}{d\mathbf{y}_k}=  \bcm \mm\mu^\top & \mm\nu^\top \ecm   \bcm
    \frac{d \mathbf{f}^*}{d  \mathbf{y}_k } \\  \frac{d \mathbf{g}^*}{d  \mathbf{y}_k }     \ecm \\
    &= \frac{1}{ \frac{1}{M}+\frac{1}{N} }\bcm \frac{1}{M}\mathbb{1}^\top_M & \frac{1}{N} \mathbb{1}^\top_N \ecm   \bcm 
    \mathbf{e}_k  (\mathbf{B}_k\mathbb{1}_N)^\top \\
    \mathbf{B}_k^\top\ecm  = \sum_{j=1}^N \frac{\partial C_{kj}}{\partial \mathbf{y}_k}\Pi^*_{kj}.
\end{align*}
If we define the third-order tensor $\mathcal{B}$ by stacking the matrices $\{\mathbf{B}_k\}_{k=1}^M$, that is  the $k$-th component $\mathcal{B}_{k\cdot \cdot}= \mathbf{B}_k$, then the linear system \eqref{linear-sys-k} can be further vectorized as 
\begin{align}  \label{linear-sys}
\mathbf{H}^* \bcm
    \frac{d \mathbf{f}^*}{d \mathbf{Y} } \\  \frac{d \mathbf{g}^*}{d \mathbf{Y} } 
    \ecm = \mathcal{R}, \ \   \text{with}\quad  \mathcal{R}: = \bcm \text{diag}(\mathcal{B} \cdot\mathbb{1}_N) \\ \mathcal{B}^{\top} \ecm  \in \mathbb{R}^{(M+N)\times M\times d}
\end{align}
where $\text{diag}(\mathcal{B} \cdot\mathbb{1}_N)$ is a third-order tensor, with  $\text{diag}(\mathcal{B} \cdot\mathbb{1}_N)_{kks} = (\mathcal{B}\cdot \mathbb{1}_N)_{ks} $ and zeros for the other entries, and  $\mathcal{B}^{\top}$ is the transpose of permutation of indices defined by $\left(\mathcal{B}^{\top}\right)_{ijk}=\mathcal{B}_{kij}$. 
Hence, by plugging \eqref{linear-sys} back to \eqref{deriv_AllSample}, the derivative of $\text{OT}_{\epsilon}(\mathbf{C}, \mm\mu, \mm\nu)$ with respect to the source data $\mathbf{Y}$ in the tensor form is equal to
$\frac{d\mathrm{OT}_\epsilon (\mathbf{C}, \mm\mu, \mm\nu)}{d\mathbf{Y}}=  \mathcal{B}\cdot \mathbb{1}_N$. 
\end{proof}

 This theorem not only implies that automatic differentiation is unnecessary for computing the gradient of EOT distance but also that solving the large linear system \eqref{linear-sys} for $\frac{d \mathbf{f}^*}{d \mathbf{y}_k }$ and $\frac{d \mathbf{g}^*}{d \mathbf{y}_k }$ is not needed. The analytical expression of the gradient follows directly from the cost matrix $\mathbf{C}_\theta$ and optimal coupling matrix $\mm\Pi^*$ from the Sinkhorn algorithm, which can significantly speed up the computation. This theorem provides a different approach to derive the same analytical result other than the Danskin's theorem. 
 We emphasize that the gradient formula  \eqref{derivate_data_Y} is generic and applies to any choice of cost matrix $\mathbf{C}$ without requiring $\mm\mu$ and $\mm\nu$ to be uniform. In particular, if the Sinkhorn iteration stops early and the coupling matrix $\widehat{\mm\Pi}$ is suboptimal, 
 the formula is still exact for $ \frac{d\mathrm{OT}_\epsilon (\mathbf{C}, \tilde{\mm \mu}, \tilde {\mm\nu})}{d\mathbf{Y}}$ given suboptimal coupling matrix $\widehat{\mm\Pi}$ and associated marginals $\tilde{\mm\mu}$ and $\tilde{\mm\nu}$, where  $\sum_{j}\widehat {\Pi}_{ij}=\widetilde{\mu}_i$ and $\sum_{i}\widehat {\Pi}_{ij}=\widetilde{\nu}_j$.

\bigskip
\noindent
\textbf{Gradient of $\widetilde{\text{OT}}_\epsilon (\mathbf{C}, \mm\mu,  \mm\nu)$. } {\color{black}{Unlike Danskin's approach, the above techniques also applies for computing the gradient of the Sinkhorn distance $\widetilde{\text{OT}}_\epsilon (\mathbf{C},  \mm\mu,  \mm\nu)$.}} Although it does not yield an analytical expression, it helps significantly reduce the computational cost by avoiding directly solving the linear system \eqref{linear-sys}, which is costly since the right-hand-side is a large third-order tensor $\mathcal{R}$. 
{\color{black}{Specifically, we compute the derivative with respect to the source data $\mathbf{Y}$, whereas prior works, such as \cite{luise2018differential}, have focused on computing the derivative with respect to the measures $\mm{\mu}$ or $\mm{\nu}$. This distinction represents the following  \underline{novel result}. }}

\begin{proposition}\label{prop:Sink_dist}
The gradient of Sinkhorn distance in \eqref{sink_distance} with respect to the source data $\mathbf{Y}$ is 
\begin{align}\label{sink_dist_derivative}
     &\frac{d \widetilde{\mathrm{OT}}_\epsilon (\mathbf{C},  \mm\mu,  \mm\nu)}{d \mathbf{Y}}  = 
     \left(\mathcal{B}- \frac{\mathcal{A}}{\epsilon}\right)\cdot \mathbb{1}_N 
    +  \frac{1}{\epsilon}  \mm r^\top  \mathcal{R},
\end{align}
where $\mathcal{B}$ and $\mathcal{R}$ are defined in \eqref{linear-sys}, the third-order tensor $\mathcal{A}$ is $\mathcal{A}_{ksj}=\frac{\partial C_{kj}}{\partial (\mathbf{y}_{k})_s}C_{kj}\Pi^*_{kj}$. The vector $\mm r$ is the solution of the linear system, 
$    \mathbf{H}^*\mm r =\bcm  \mm a \\  \mm b \ecm
$, 
where 
vectors $\mm a$ and $\mm b$ have entries
$a_i = \sum_{j=1}^{M}  C_{ij}\Pi^*_{ij}$ and  $ b_j = \sum_{i=1}^{N}  C_{ij}\Pi^*_{ij}$.
\end{proposition}
\begin{proof}
From \eqref{sink_distance} and \eqref{EOT_Pi_derivative}, the derivative of Sinkhorn distance with respect to the data $\mathbf{y}_k$ is 
\begin{align*}
     \frac{d \widetilde{\mathrm{OT}}_\epsilon (\mathbf{C},  \mm \mu,  \mm \nu)}{d \mathbf{y}_k} &=  \sum_{j=1}^N \left(1-\frac{C_{kj}}{\epsilon}\right)\Pi^*_{kj}\frac{\partial C_{kj}}{\partial \mathbf{y}_k} +  \frac{1}{\epsilon}\left(\sum_{i=1}^M a_i\frac{\partial f_i^*}{\partial \mathbf{y}_k} + \sum_{j=1}^N b_j\frac{\partial g_j^*}{\partial \mathbf{y}_k}\right),
\end{align*}
where $a_i = \sum_{j=1}^{N} C_{ij}\Pi^*_{ij}$, and  $b_j = \sum_{i=1}^{M}  C_{ij}\Pi^*_{ij}$.
The main computation burden is the second term.  To compute $\left(\sum_{i} a_i\frac{\partial f_i^*}{\partial \mathbf{y}_k} + \sum_j b_j\frac{\partial g_j^*}{\partial \mathbf{y}_k}\right)$, a common practice is 
to solve the linear system \eqref{linear-sys} to get all $\frac{\partial f^*}{\partial \mathbf{y}_k}$ and $\frac{\partial g^*}{\partial \mathbf{y}_k}$, which is computationally expensive because the right-hand-side of \eqref{linear-sys} is  a third-order tensor $\mathcal{R}$. 

\noindent
Although we do not have \eqref{key-observation} as $\mm a$ and $\mm b$ are no longer the marginal constraints of $\mm\Pi^*$, we can yet apply the same technique to reduce the computation costs by solving the following matrix-vector form of linear system only once
 $ \mathbf{H}^*   \mm r  = \bcm \mm a \\ \mm b \ecm \in \mathbb{R}^{(M+N)}$. 
That is, we only need to solve for one time the matrix-vector form of the linear system, for the column vector $  \mm r$. 
This can be efficiently solved by conjugate gradient method with early stopping.

Notice that $\bcm \mathbb{1}^\top_{M} &-\mathbb{1}^\top_N\ecm \bcm \mm{a} \\ \mm{b}\ecm  =0 $, and  $\text{span}\left(\bcm \mathbb{1}_{M}\\-\mathbb{1}_N\ecm\right)$ is the kernel  of $\mathbf{H}^*$ as proved in Theorem~\ref{prop:H_general_Pi},
so $\bcm \mm a \\ \mm b \ecm$ is in the span of $\mathbf{H}^*$,
therefore the linear system above always has a solution. 

\noindent
Notice that 
\begin{align*}
\sum_{i} a_i\frac{\partial f_i^*}{\partial \mathbf{Y}} + \sum_j b_j\frac{\partial g_j^*}{\partial \mathbf{Y}}
=\bcm \mm a^\top &\mm b^\top \ecm \bcm\frac{d \mathbf{f}^*}{d\mathbf{Y}}\\ \frac{d \mathbf{g}^*}{d\mathbf{Y}}\ecm 
&=  \mm r^\top \mathbf{H}^* \bcm\frac{d \mathbf{f}^*}{d\mathbf{Y}}\\ \frac{d \mathbf{g}^*}{d\mathbf{Y}}\ecm= \mm r^\top   \mathcal{R}.
\end{align*}
Hence, we prove \eqref{sink_dist_derivative}.
\end{proof}

\subsection{Computation of the Hessian}\label{subsec:Hessian_EOT}
In this subsection, 
we \textit{analytically} compute the Hessian of the loss function with respect to the source data $\bY$.  We mainly focus on the EOT distance.

\begin{theorem}    \label{theo:Hessian}
The second-order derivative of EOT distance $\mathrm{OT}_{\epsilon}(\mathbf{C},  \mm \mu,  \mm \nu)$ with respect to the source data $\mathbf{Y}$ is given by the fourth order tensor $\mathcal{T}\in \R^{M\times d\times M\times d}$
\begin{equation}\label{eq:hessian_OTe}
\mathcal{T}_{ktsl}  =
    \frac{1}{\epsilon} \sum_{i,j=1}^{M+N}\mathcal{R}_{ikt}  (\mathbf{H}^*)^\dag_{ij} \mathcal{R}_{jsl} +\mathcal{E}_{ktsl}, \text{ for } k,s=1,\dots,M \text{ and } t,l=1,\dots, d, 
   \end{equation}
   where $(\mathbf{H}^*)^\dag$ is the Moore–Penrose inverse of matrix $ \mathbf{H}^* \in \R^{(M+N)\times (M+N)}$ defined in \eqref{eq:H}, and $\mathcal{R}= \bcm \text{diag}(\mathcal{B} \cdot\mathbb{1}_N) \\ \mathcal{B}^{\top } \ecm  \in  \R^{(M+N)\times M \times d}$ is the right-hand-side third-order tensor defined in \eqref{linear-sys}. The fourth-order tensor $\mathcal{E}\in \R^{M\times d\times M\times d}$ is defined as 
   \begin{align} \label{eq:Hessian_E}
    \mathcal{E}_{kt s l } =
    \begin{cases}
    \sum_{j=1}^N\Pi_{kj}^* \left(\left(\frac{\partial^2 C_{kj}}{\partial \mathbf{y}_k^2}\right)_{tl}- \frac{1}{\epsilon}\frac{\partial C_{kj}}{\partial (\mathbf{y}_k)_t}\cdot\frac{\partial C_{kj}}{\partial (\mathbf{y}_k)_l}\right),   & \text{if } k= s \\
    0, & \text{Otherwise}
    \end{cases} 
\end{align}
for $k,s = 1,\dots, M$ and $t, l=1, \dots, d$.
    \end{theorem}
\begin{proof}
    With the help of the gradient given by \eqref{derivate_data}, 
    we know 
\begin{equation}\label{Hessian_sample}
    \frac{d^2\mathrm{OT}_{\epsilon}(\mathbf{C},  \mm\mu,  \mm\nu)}{d\mathbf{y}_s d\mathbf{y}_k} =   \sum_{j=1}^N  \left(\frac{\partial C_{kj}}{\partial \mathbf{y}_k}\right)^\top \frac{\partial\Pi_{kj}^*}{\partial \mathbf{y}_s} + \delta_{ks} \sum_{j=1}^N  \frac{\partial^2 C_{kj}}{\partial \mathbf{y}_k^2}\Pi_{kj}^*. 
\end{equation}
For $k\neq s$, we have $\delta_{ks}=0$, so the term is expanded as
\begin{align*}
    &\sum_{j=1}^N  \left(\frac{\partial C_{kj}}{\partial  \mathbf{y}_k}\right)^\top\frac{\partial\Pi_{kj}^*}{\partial \mathbf{y}_s} = \sum_{j=1}^{N} \frac{\Pi^*_{kj}}{\epsilon}\left( \frac{\partial C_{kj}}{\partial \mathbf{y}_k}\right)^\top  \left( \frac{\partial f^*_k}{\partial \mathbf{y}_s}+\frac{\partial g^*_j}{\partial \mathbf{y}_s}\right)  \\
    & =  \bcm 
    \mm{e}_k  (\mathbf{B}_k\mathbb{1}_N)^\top \\
    \mathbf{B}_k^\top\ecm^\top \frac{ (\mathbf{H}^*)^\dag}{\epsilon}  \bcm 
    \mm{e}_s  (\mathbf{B}_s\mathbb{1}_N)^\top \\
    \mathbf{B}_s^\top\ecm,
    \end{align*}
where we use the fact that $ \bcm \frac{ \partial \mathbf{f}^*}{\partial \mathbf{y}_s} \\ \frac{\partial \mathbf{g}^*}{\partial \mathbf{y}_s}\ecm$ is the solution of the linear system \eqref{linear-sys-k}, so it can be expressed as  
  $  \bcm \frac{ \partial \mathbf{f}^*}{\partial \mathbf{y}_s} \\ \frac{\partial \mathbf{g}^*}{\partial \mathbf{y}_s}\ecm = (\mathbf{H}^*)^\dag \bcm 
    \mathbf{e}_s  (\mathbf{B}_s\mathbb{1}_N)^\top \\
    \mathbf{B}_s^\top\ecm$.
For $k=s$, the term is expanded as follows
\begin{align*}
&\sum_{j=1}^N \frac{\partial^2 C_{kj}}{\partial \mathbf{y}_k^2}\Pi_{kj}^* +  \sum_{j=1}^N  \left(\frac{\partial C_{kj}}{\partial \mathbf{y}_k}\right)^\top\frac{\partial\Pi_{kj}^*}{\partial \mathbf{y}_k}\\ 
&\,= \sum_{j=1}^N 
\frac{\Pi^*_{kj}}{\epsilon} \left(\frac{\partial C_{kj}}{\partial \mathbf{y}_k}\right)^\top \left( \frac{\partial f^*_k}{\partial \mathbf{y}_k}+\frac{\partial g^*_j}{\partial \mathbf{y}_k}\right)+ \sum_{j=1}^N\Pi_{kj}^*\left(\frac{\partial^2 C_{kj}}{\partial \mathbf{y}_k^2} - \frac{1}{\epsilon} \left(\frac{\partial C_{kj}}{\partial \mathbf{y}_k}\right)^\top \left(\frac{\partial C_{kj}}{\partial \mathbf{y}_k}\right)\right) \\
&\,= \bcm 
    \mm{e}_k  (\mathbf{B}_k\mathbb{1}_N)^\top \\
    \mathbf{B}_k^\top\ecm^\top \frac{(\mathbf{H}^*)^\dag}{\epsilon} \bcm 
    \mm{e}_k  (\mathbf{B}_k\mathbb{1}_N)^\top \\
    \mathbf{B}_k^\top\ecm  +\sum_{j=1}^N\Pi_{kj}^*\left(\frac{\partial^2 C_{kj}}{\partial \mathbf{y}_k^2} -  \frac{(\frac{\partial C_{kj}}{\partial \mathbf{y}_k})^\top (\frac{\partial C_{kj}}{\partial \mathbf{y}_k})}{\epsilon}\right).
\end{align*}
In vector form, we obtain the Hessian in \eqref{eq:hessian_OTe}. 
\end{proof}

   Given the cost matrix $\mathbf{C}$ and the optimal coupling matrix $\mm\Pi^*$ from Sinkhorn algorithm, the Hession tensor is analytically calculated once the large linear system \eqref{linear-sys} is numerically solved. Unlike the previous approaches, the solution of the linear system is no longer needed to be differentiated in order to obtain the second-order derivatives, 
   so our direct analytical Hessian expression could significantly speed up the computation with less memory burden. 
   Similar to the case of the first-order derivative, the analytical expression for the Hessian is generic and applicable to any choice of cost function. It does not require $\mm \mu$ and $\mm \nu$ to be uniform either.
   If the Sinkhorn iterations stops early with $\widehat{\mm\Pi}$ being suboptimal, then the Hessian expression \eqref{eq:hessian_OTe} is still exact for $\frac{d^2\mathrm{OT}_{\epsilon}(\mathbf{C}, {\tilde{\mm\mu}}, {\tilde{\mm\nu}})}{d\mathbf{Y}^2}$ with suboptimal coupling matrix $\widehat{\mm\Pi}$ and associated marginals $\tilde{\mm\mu}$ and $\tilde{\mm\nu}$, where $\sum_{j}\widehat {\Pi}_{ij}=\widetilde{\mu}_i$ and $\sum_{i}\widehat {\Pi}_{ij}=\widetilde{\nu}_j$. 

With the analytical expression of Hessian, the following result show that the sum of the first index of the Hessian tensor is only dependent on the marginal probability vector $\mm\mu$. Later, it will be used as a marginal error to verify the accuracy of Hessian in the numerical implementation.

\begin{proposition} \label{proposition:Hessian_sum}
    If $C_{kj}=\|\mathbf{y}_k-\mathbf{y}_j^*\|_2^2$ for each $k,j$, 
    then 
        $\sum_{k=1}^M    \mathcal{T}_{k\cdot s\cdot } = 2\mu_s \mathbb{I}_d$,
{\color{black}{where $\mathbb{I}_d$ is the identity matrix of size $d$.}}
\end{proposition}
\begin{proof}
   Since the cost is square distance $C_{kj}=\|\mathbf{y}_k-\mathbf{y}_j^*\|_2^2$, then $\frac{\partial C_{kj}}{\partial \mathbf{y}_k} = 2(\mathbf{y}_k-\mathbf{y}_j^*)$ and $\frac{\partial ^2 C_{kj}}{\partial  \mathbf{y}_k^2}=2\mathbb{I}_{d}$. 
   The first term in the second-order derivative \eqref{Hessian_sample} expands as 
   \begin{align*}
       \sum_{k=1}^M \sum_{j=1}^N  \left(\frac{\partial C_{kj}}{\partial \mathbf{y}_k}\right)^\top \frac{\partial\Pi_{kj}^*}{\partial \mathbf{y}_s} & = 2\sum_{k=1}^M \mathbf{y}_k^\top \sum_{j=1}^N\frac{\partial\Pi_{kj}^*}{\partial \mathbf{y}_s} - 2\sum_{j=1}^N (\mathbf{y}_j^*)^\top  \sum_{k=1}^M\frac{\partial\Pi_{kj}^*}{\partial \mathbf{y}_s} =\mathbb{0}_{d\times d}. 
   \end{align*}
   The last equal sign is because that  $\sum_j \Pi^*_{kj}=\mu_k$ and  $\sum_k \Pi^*_{kj}=\nu_j$, as well as  $\sum_{j}\left(\frac{\partial \Pi^*_{kj}}{\partial \mathbf{y}_s}\right)=\frac{\partial \sum_{j} \Pi^*_{kj}}{\partial \mathbf{y}_s}
   $ and $\sum_{k}\left(\frac{\partial \Pi^*_{kj}}{\partial \mathbf{y}_s}\right)=\frac{\partial \sum_{k} \Pi^*_{kj}}{\partial \mathbf{y}_s}
   $.  
   
The remaining term is 
       $\sum_{j=1}^N \frac{\partial^2 C_{sj}}{\partial \mathbf{y}_s^2}\Pi^*_{sj} = 2\mathbb{I}_{d}\sum_{j=1}^N \Pi^*_{sj} = 2\mu_s\mathbb{I}_{d}.$
\end{proof}

\subsection{Solve the linear systems with truncated SVD}\label{subsec:TrunSVD}

A major challenge in computing the gradient of Sinkhorn distance in \eqref{sink_dist_derivative} and the Hessians of {\color{black}{Entropy-regularized}} distance in \eqref{eq:hessian_OTe} is that pseudo-inverse $(\mathbf{H}^*)^\dag$ can severely amplify the rounding errors or early stopping errors when the matrix $\mathbf{H}^*$ is ill-conditioned. The pseudo-inverse comes from the solutions to the linear systems \eqref{linear-sys-k} and \eqref{linear-sys} for computation of first- and second-order derivatives. Thus, instead of directly using $(\mathbf{H}^*)^\dag$, it is important to study the spectrum of the matrix $\mathbf{H}^*$ and regularize properly when the linear systems are ill-posed.

Analytical and empirical results in the next section show that the $\mathbf{H}$-matrix is often ill-conditioned when the Sinkhorn regularization parameter $\epsilon$ is small or when the optimal coupling matrix $\mm\Pi^*$ is close to a permutation. Such ill-conditioned $\mathbf{H}$-matrices result in numerically unstable solutions when $\mm\Pi^*$ is slightly perturbed to $\widehat{\mm\Pi}$ due to the early stopping of the Sinkhorn iterations. This is also the exact reason that the previous implicit differentiation approach failed due to numerical instability.

We tackle the potential ill-posedness by truncated Singular value decomposition (TSVD). We truncate the H-matrix's spectrum up to the $K$-th largest eigenvalue, with $K=\max\{j: \frac{\lambda_j}{\lambda_1}>\alpha\}$, where $\{\lambda_j\}_{j=1}^{M+N}$ are the eigenvalues of $\mathbf{H}^*$ in descending order. 
 In practice, we use the LAPACK's DGELSD algorithm building in the least-square solver \cite{anderson1999lapack} and set $\alpha = 10^{-10}$. The algorithm for gradient and Hessian computation is summarized in Algorithm~\ref{Alg_Hess}.

\begin{algorithm}
\caption{ Computation of gradient and Hessian of EOT distance $\text{OT}_{\epsilon}(C, \mm{\mu},\mm{\nu})$ with respect to the source data $\mm Y$.
\label{Alg_Hess}}
\textbf{Input:} Optimized $\Pi^*$, cost matrix $C$, entropic regularization strength $\epsilon$, source data $\mm Y$ and singular value threshold $\alpha$.\\
\textbf{Output:} Gradient with respect to $\mm Y$: $\frac{d\text{OT}_\epsilon(C, \mm \mu, \mm \nu)}{d \mm Y}\in \mathbb{R}^{M\times d}$, 
Hessian with respect to $\mm Y$: $\mathcal{T}\in \mathbb{R}^{M \times d\times M\times d}$  
\begin{algorithmic}[1]
\State Compute the marginal probability vector $\mm \mu\in \mathbb{R}^M$ and $\mm \nu\in \mathbb{R}^N$: $\mm \mu \leftarrow \Pi^* \mathbb{1}_N$ and $\mm \nu \leftarrow (\Pi^*)^\top \mathbb{1}_M $.

\State Compute matrix $H\in \mathbb{R}^{(M+N)\times (M+N)}$: $H \leftarrow \bcm \text{diag}(\mm \mu) & \Pi^* \\ (\Pi^*)^\top & \text{diag}(\mm \nu)\ecm$.

\State Compute third-order tensor $\mathcal{B}\in \mathbb{R}^{M\times d\times N}$: $\mathcal{B}_{ksj}\leftarrow  \frac{\partial C_{kj}}{\partial (\mm{y}_k)_s}\Pi^*_{kj}$,  compute third order tensor $\mathcal{R}\in \mathbb{R}^{(M+N)\times M \times d}$: $\mathcal{R}\leftarrow \bcm \text{diag}(\mathcal{B}\cdot \mathbb{1}_N) \\ \mathcal{B}^\top\ecm$ and compute the fourth-order tensor $\mathcal{E}$ in \eqref{eq:Hessian_E}.

\State Compute the gradient: $\frac{d\text{OT}_\epsilon(C, \mm \mu, \mm \nu)}{d \mm Y} = \mathcal{B}\cdot \mathbb{1}_N$.

\State Compute the truncated singular value decomposition of $H$-matrix up to the $K$-th largest eigenvalue: $H\approx U_K\Lambda_K U_K^\top$ with $K=\max\{j: \frac{\lambda_1}{\lambda_j}>\alpha\}$.

\State Approximate $\bcm \frac{d \mm f^*}{d \mm Y} \\ \frac{d \mm g^*}{d \mm Y} \ecm$: $\bcm \frac{d \mm f^*}{d \mm Y} \\ \frac{d \mm g^*}{d \mm Y} \ecm \leftarrow U_K \Lambda_K^{-1}U_K^\top \mathcal{R}$.
\State  Compute the Hessian $\mathcal{T}$: $\mathcal{T}_{ktsl}\leftarrow \frac{1}{\epsilon}\sum_{i=1}^{M+N} \mathcal{R}_{ikt}\bcm \frac{d \mm f^*}{d \mm Y} \\ \frac{d \mm g^*}{d \mm Y} \ecm_{isl}+\mathcal{E}_{ktsl}$.
\end{algorithmic}
\end{algorithm}

\section{Spectral analysis of the H-matrix}\label{sec:spectrum}

{\color{black}{
This section analyzes the spectrum of the matrix $\mathbf{H}$ in \eqref{eq:H} and addresses the potential ill-posedness in solving the linear systems in \eqref{linear-sys}, which are critical for computing the gradient and Hessian. 

Recall that $\mathbf{H}$ is a block matrix consisting of the coupling matrix $\mm\Pi \in \mathbb{R}^{M \times N}$ and diagonal matrices of its two marginal vectors $\mm\mu$ and $\mm\nu$. We assume that the coupling matrix $\mm\Pi$ is known, either estimated via the Sinkhorn algorithm or computed analytically. We focus on two common types of coupling matrices: (i) a positive coupling matrix, where $\Pi_{ij} > 0$ for all $i, j$, and (ii) a coupling matrix with uniform marginal distributions, where $\mm\mu = \mathbb{1}_M / M$ and $\mm\nu = \mathbb{1}_N / N$. The first type typically arises from the Sinkhorn algorithm, as all entries of the optimal coupling matrix are positive, as shown in \eqref{KOT_dual} \cite{feydy2020geometric,peyre2019computational}. The second type is common in applications with randomly sampled data. In entropic OT problems, the second type can be viewed as a special case of the first.

Our main contribution here is to establish upper and lower bounds for the condition number of $\mathbf{H}$ based on the spectral gap of $\mm\Pi^\top \mm\Pi$. We also provide examples to illustrate that $\mathbf{H}^*$ can be ill-conditioned. So the truncated SVD regularization, as discussed in Section \ref{subsec:TrunSVD}, is necessary to mitigate numerical instability.

First, in Section \ref{sec:general-H}, we prove in Theorem \ref{prop:H_general_Pi} that $\mathbf{H}$ is singular, with a simple zero eigenvalue, for any positive coupling matrix. Although removing one of the rows can eliminate the singularity, the linear system may still remain ill-posed.

Second, in Section \ref{sec:H-unif_Pi}, we analytically compute the eigenvalues of $\mathbf{H}$ for coupling matrices with uniform marginal distributions (including permutation matrices). In particular, we derive lower and upper bounds for the condition number of $\mathbf{H}$ in terms of the spectral gap of $\mm\Pi^\top \mm\Pi$, as established in Theorem \ref{thm:condi_H1}. These bounds are extended to the $\mathbf{H}$-matrix computed by the Sinkhorn algorithm with early stopping in the entropic optimal transport (EOT) distance.

Finally, in Section \ref{subsec:NumericalSpecStudy}, we derive the sharp asymptotic behavior of the condition number for one example in terms of the regularization strength $\epsilon$ and the sample size $N$, demonstrating that this linear system can be highly ill-conditioned. Interestingly, we observe similar asymptotic behavior empirically when applying to other randomly sampled datasets.
}
}
     
\subsection{General positive coupling matrices}
\label{sec:general-H}
We show first that the  $\mathbf{H}$ is singular with a simple zero eigenvalue for any positive coupling matrix $\mm\Pi$. As a result, we call 
$\kappa(\mathbf{H})=\frac{\lambda_1(\mathbf{H})}{\lambda_{N+M-1}(\mathbf{H})}$, 
 the \emph{condition number} of the matrix $\mathbf{H}$. 
 
\begin{theorem}[\textbf{Simple zero eigenvalue for the $\mathbf{H}$-matrix}]
\label{prop:H_general_Pi}
For any positive coupling matrix $\mm\Pi \in \R^{M\times N}_{>0}$,  
the smallest eigenvalue of $\mathbf{H}$
is zero and it is simple, with eigenvector $\mathbf{q}_0=\bcm \mathbb{1}_M \\ -\mathbb{1}_N \ecm$. 
\end{theorem}

Its proof relies on the next Perron-Frobenius type lemma, which shows that the largest eigenvalue of the matrix $\text{diag}(\mm\nu)^{-1} \mm\Pi^\top \, \text{diag}(\mm\mu)^{-1} \mm\Pi$ is $1$ and is simple. We postpone the proof of the lemma to supplementary materials. 

\begin{lemma}\label{lemma:normalize_Pi}
For any coupling matrix with strictly positive entries $\mm\Pi\in \R^{M\times N}_{>0}$, the largest eigenvalue of matrix $\text{diag}(\mm\nu)^{-1} \mm\Pi^\top \, \text{diag}(\mm\mu)^{-1} \mm\Pi$ is $\lambda=1$ and has multiplicity one, with eigenvector $ \mathbb{1}_N$. Similarly, the largest eigenvalue of the matrix $\text{diag}(\mm\mu)^{-1} \mm\Pi \, \text{diag}(\mm\nu)^{-1} \mm\Pi^\top$ is $\lambda=1$ and has multiplicity one, with eigenvector $ \mathbb{1}_M$.
\end{lemma}

\begin{proof}[\textbf{Proof of Theorem \ref{prop:H_general_Pi}}]
It is clear that $0$ is an eigenvalue with eigenvector $\mathbf{q}_0$. 
    We show next that any other eigenvector of $0$ must be  $\mathbf{q}_0$ up to a scalar factor.  
    Note that the vector $\bcm \mm w \\ \mm v \ecm$ is the eigenvector of $0$ if and only if
        \begin{align} \label{eq:Heigen0_eq}
           \bcm\text{diag}( \mm \mu) & \mm\Pi \\ \mm\Pi^\top & \text{diag}(\mm \nu)
    \ecm \bcm \mm w \\ \mm v \ecm = \mathbb{0}_{M+N} \Leftrightarrow   \bcm  \text{diag}( \mm\mu) \mm w + \mm\Pi \mm v = 0  \\  \text{diag}( \mm \nu) \mm v + \mm\Pi^\top \mm w  =0 \ecm, 
    \end{align}
which is equivalent to 
\begin{equation*}
     \begin{aligned}
     \mm v&= - \text{diag}(\mm \nu)^{-1} \mm\Pi^\top \mm w=  \text{diag}(\mm \nu)^{-1} \mm\Pi^\top \text{diag}(\mm \mu)^{-1} \mm\Pi \mm v, \\
     \mm w&= - \text{diag}(\mm \mu)^{-1} \mm\Pi \mm v=  \text{diag}(\mm \mu)^{-1} \mm\Pi \text{diag}(\mm \nu)^{-1} \mm\Pi^\top \mm w.
       \end{aligned}
       \end{equation*}    
     By Lemma \ref{lemma:normalize_Pi}, $\mm v= a  \mathbb{1}_N$ and $\mm w=  b \mathbb{1}_M$ for some nonzero constant $a,b$ for each equation to hold. 
     Plugging back to \eqref{eq:Heigen0_eq}, we have $a=-b$, and $\bcm \mm w \\ \mm v \ecm = a \mathbf{q}_0 $. 
     Thus, the zero eigenvalue of $\mathbf{H}$ is simple. 
\end{proof}

A quick corollary of the above lemma is that the smallest eigenvalue of $\mathbb{I}-\text{diag}(\mm \nu)^{-1} \mm\Pi^\top \text{diag}(\mm \mu)^{-1} \mm\Pi $ is zero and is simple. As a result, we obtain an invertible matrix after dropping one of the rows. 
\begin{corollary}\label{cor:drop1row}
For any coupling matrix $\mm\Pi$, $\text{diag}(\mm{\overline \nu}) - \overline {\mm\Pi}^\top \text{diag}(\mm\mu)^{-1} \mm\Pi $ is invertible, where $\mm{\overline \nu}$ and $\overline{\mm\Pi}$ are the arrays after dropping the last rows.
\end{corollary}
The invertibility of the above matrix has been used in \cite{xie2021hypergradient, luise2018differential, blondel2022efficient, eisenberger2022unified} to remove the zero eigenvalues of $\mathbf{H}^*$ in the computation of the gradient of Sinkhorn distance \eqref{sinkhorn_dist}. However, the ill-posedness in solving the linear systems is not addressed as the other eigenvalues are still close to $0$. Next, we analyze the spectrum of the $\mathbf{H}$-matrix, which guides the treatment of the ill-posedness.

\subsection{Coupling matrices with uniform marginal distributions}\label{sec:H-unif_Pi}
Coupling matrices with uniform marginal distributions are of particular interest as they arise in many applications.
In this section, we compute the eigenvalues of the $\mathbf{H}$-matrices for coupling matrices $\mm\Pi$ with uniform marginal distributions (including the permutation matrices), {\color{black}{i.e., $\mm\Pi\mathbb{1} = \frac{1}{M}\mathbb{1}_M$ and $\mm\Pi^\top\mathbb{1} = \frac{1}{N}\mathbb{1}_N$.}} These eigenvalues shed light on the root of the ill-conditionedness of the $\mathbf{H}$-matrices. In particular, if the coupling matrix is entrywise positive, 
we provide upper and lower bounds for the condition number of the $\mathbf{H}$-matrix in terms of the spectral gap of $\mm\Pi^\top \mm\Pi$ in Theorem \ref{thm:condi_H1}.

\begin{proposition}[\textbf{Eigenvalues of $\mathbf{H}$ and singular values of $\mm\Pi$}] 
\label{H-Pi_eig_svd}	
Let $\mm\Pi\in [0,1]^{M\times N}$ be a (not necessarily positive) coupling matrix with uniform marginal distributions. Let $M\leq N$ and assume $\mm\Pi$ has rank $M$. Then, the eigenvalues (in descending order) and eigenvectors of $\mathbf{H}$ defined in \eqref{eq:H} are, for $j=1, \dots, M $, 
\begin{equation}\label{eq:eigH}
\begin{aligned}
    \lambda_{j}(\mathbf{H}) &= \frac{\left(\frac{M+N}{MN}\right)+ \sqrt{\left(\frac{M-N}{MN}\right)^2 +4\sigma_j(\mm\Pi)^2}}{2}, \quad \mathbf{q}_j  =\bcm  \frac{\kappa_j}{\sqrt{1+\kappa_j^2}} \mathbf{u}_j  \\ \frac{1}{\sqrt{1+\kappa_j^2}}\mathbf{v}_j \ecm, \\
    \lambda_{M+1}(\mathbf{H})&=\dots = \lambda_{N}=\frac{1}{N}, \quad \mathbf{q}_{i} = \bcm \mathbb{0} \\ \mathbf{v}_i \ecm, \, M<i\leq N, \\
    \lambda_{N+M+1-j}(\mathbf{H}) &= \frac{\left(\frac{M+N}{MN}\right)- \sqrt{\left(\frac{M-N}{MN}\right)^2 +4\sigma_j^2(\mm\Pi)}}{2}, \quad\mathbf{q}_{N+M+1-j} = \bcm  \frac{-1}{\sqrt{1+\kappa_j^2}}\mathbf{u}_j  \\ \frac{\kappa_j}{\sqrt{1+\kappa_j^2}} \mathbf{v}_j \ecm,  
\end{aligned}
\end{equation}
where $\sigma_j(\mm\Pi)$ is the $j$-th singular value of $\mm\Pi$ (in descending order) with $\mathbf{u}_j$ and $\mathbf{v}_j$ being the right and left singular vectors, 
   $\{\mathbf{v}_i\}_{i=M+1}^N$ are the  $N - M$ orthogonal vectors in the kernel of $\mm\Pi$, and these scalars $\kappa_j =\left(\frac{\frac{1}{M}-\frac{1}{N}}{2\sigma_j(\mm\Pi)}\right)+ \sqrt{\left(\frac{\frac{1}{M}-\frac{1}{N}}{2\sigma_j(\mm\Pi)}\right)^2 +1}$.

In particular, when $M=N$ and $\mm\Pi$ is a permutation coupling matrix, i.e., its rows are a permutation of the rows of $\frac{1}{N}I_N$. Then, the eigenvalues of $H(\mm\Pi)$ are $\lambda_1(H)=\cdots = \lambda_N(H)= \frac 2 N$ and $\lambda_{N+1}(H)=\cdots= \lambda_{2N}(H)=0$. 
\end{proposition}

\begin{proof}[\textbf{Proof of Proposition \ref{H-Pi_eig_svd}}]

Denote an eigen-pair $H$ by $\left(\lambda, \bcm \mathbf{u}\\  \mathbf{v}\ecm\right)$, where $\mathbf{u} \in \mathbb{R}^{M}$ and $\mathbf{v} \in \mathbb{R}^{N}$; that is,  $
   \bcm \frac{1}{M}\mathbb{I}_M & \mm\Pi \\ (\mm\Pi)^\top & \frac{1}{N}\mathbb{I}_N \ecm  \, \bcm \mathbf{u}\\  \mathbf{v}\ecm = \lambda \bcm \mathbf{u}\\  \mathbf{v}\ecm$. 

Then, $\mm\Pi \mathbf{v}= (\lambda-\frac{1}{M}) \mathbf{u}$ and $\mm\Pi^\top \mathbf{u}= (\lambda-\frac{1}{N}) \mathbf{v}$. This implies that 
\begin{align*} 
   \left( \mm\Pi (\mm\Pi)^\top \right) \mathbf{u} = \left(\lambda - \frac{1}{M}\right)  \left(\lambda - \frac{1}{N}\right)\mathbf{u},  \quad \left( (\mm\Pi)^\top \mm\Pi\right)\mathbf{v} =\left(\lambda - \frac{1}{M}\right)  \left(\lambda - \frac{1}{N}\right) \mathbf{v}.
\end{align*}
That is, $\mathbf{u}$ and $\mathbf{v}$ are the right and left singular vectors of $\mm\Pi$, corresponding to a singular value of $\mm\Pi$ satisfying 
$\sigma_j^2=\left(\lambda - \frac{1}{M}\right)  \left(\lambda - \frac{1}{N}\right)$, for some $1\leq j \leq M$. 
Solving this equation, we obtain two eigenvalues 
 $  \lambda_{\pm} = \frac{\left(\frac{1}{M}+\frac{1}{N}\right)\pm \sqrt{\left(\frac{1}{M}-\frac{1}{N}\right)^2 +4\sigma_j^2}}{2}. 
$ 

To compute the eigenvector of $\lambda_+$, let $\mathbf{v}=\mathbf{v}_j$. Then $\mm\Pi \mathbf{v} = \sigma_j \mathbf{u}_j=(\lambda_+-\frac{1}{M})\mathbf{u}$. So $\mathbf{u} = \frac{\sigma_j}{(\lambda_+-\frac{1}{M})}\mathbf{u}_j$ (recall that $\sigma_j\neq 0$ by $\text{rank}(\mm\Pi)=m$). 
 Note that  $\frac{(\lambda_+-\frac{1}{N})}{\sigma_j}=\frac{\left(\frac{1}{M}-\frac{1}{N}\right)+ \sqrt{\left(\frac{1}{M}-\frac{1}{N}\right)^2 +4\sigma_j^2}}{2\sigma_j}$; hence,  we have
\begin{align*}
 \frac{\sigma_j}{(\lambda_+-\frac{1}{M})} = \frac{(\lambda_+-\frac{1}{N})}{\sigma_j}
 =\frac{\frac{1}{M}-\frac{1}{N}}{2\sigma_j}+ \sqrt{\left(\frac{\frac{1}{M}-\frac{1}{N}}{2\sigma_j}\right)^2 +1} = \kappa_j.
\end{align*}
So an eigenvector  for $\lambda_+$ is $\bcm  \frac{\kappa_j}{\sqrt{1+\kappa_j^2}} \mathbf{u}_j  \\ \frac{1}{\sqrt{1+\kappa_j^2}}\mathbf{v}_j \ecm$. Similarly, $\lambda_-$ has $\bcm  \frac{-1}{\sqrt{1+\kappa_j^2}}\mathbf{u}_j  \\ \frac{\kappa_j}{\sqrt{1+\kappa_j^2}} \mathbf{v}_j \ecm$.

The above $\lambda_{\pm}$ account for $2M$ eigenvalues. The other $N-M$ eigenvalues correspond to $\sigma_i(\mm\Pi)=0$, which has singular vectors $\mathbf{v}_i$, for $i=M+1, \dots, N$. Thus, setting $\mathbf{u}=\mathbb{0}$ and $\mathbf{v}=\mathbf{v}_i$, we have $\mm\Pi^\top \mathbf{u}=\mathbb{0}=(\lambda-\frac{1}{N})\mathbf{v}_i$. So $\lambda=\frac{1}{N}$. This eigenvalue has multiplicity of $N-M$, with eigenvectors $\bcm \mathbb{0} \\ \mathbf{v}_i \ecm$  for $i=M+1, \dots, N$.

At last, when $M=N$ and $\mm\Pi$ is a permutation matrix, note that $\mm\Pi^\top \mm\Pi = \frac{1}{N^2}\mathbb{I}_N$. Thus, the singular values of $\mm\Pi$ are $\sigma_j=\frac 1 N$ with multiplicity $N$. Applying \eqref{eq:eigH}, we obtain that the eigenvalues of $H$ are $0$ and $\frac{2}{N}$, each with multiplicity $N$. 
\end{proof}

\begin{theorem}[\textbf{Condition number of $\mathbf{H}$}] 
\label{thm:condi_H1}
Let $\mm\Pi$ be a positive coupling matrix with uniform marginal distributions and singular values $\{\sigma_k\}$ in descending order. Then the condition number $\kappa(\mathbf{H})$ of $\mathbf{H}$ in \eqref{eq:H} is bounded by 
\begin{equation}\label{eq:condiH}
\frac{(M+N)^2}{2 M^2N^2 \rho(\mm\Pi^\top \mm\Pi) } \leq \kappa(\mathbf{H}) \leq \frac{(M+N)^2}{ M^2N^2 \rho(\mm\Pi^\top \mm\Pi) }, 
\end{equation} 
{\color{black}{where $\rho(\mm\Pi^\top \mm\Pi)$ is spectral gap of the matrix $\mm\Pi^\top \mm\Pi$.}} 
\end{theorem}

\begin{proof}[\textbf{Proof of Theorem \ref{thm:condi_H1}}]

The largest singular value of $\mm\Pi$ is $\sigma_1= \frac{1}{\sqrt{NM}}$ and it is simple. Hence, by Prop.\ref{H-Pi_eig_svd}, the largest and smallest eigenvalues of $\mathbf{H}$ are 
$
\lambda_{1}(\mathbf{H})= \frac{1}{M}+\frac{1}{N}
$, $\lambda_{M+N}(\mathbf{H})=0$, and both are simple.

\noindent
The second smallest eigenvalue of $\mathbf{H}$, denoted by $\lambda_{M+N-1}(\mathbf{H})$, can be obtained from 
$\left(\lambda - \frac{1}{M}\right)  \left(\lambda- \frac{1}{N}\right) = \sigma_2^2$, 
i.e., $\lambda^2 -\lambda_1 \lambda + \sigma_1^2- \sigma_2^2 =0 $. With $\Delta := \frac{(\sigma^2_1-\sigma^2_2)}{\lambda_1^2}$, this 
gives 
\begin{align}\label{LowestTwo_eigH}
\lambda_{M+N-1}(\mathbf{H}) &= \frac{1}{2}\lambda_1\left[1- \sqrt{1- 4\Delta }\right] = \frac{1}{2}\lambda_1\frac{4\Delta }{1+ \sqrt{1- 4\Delta}}. 
\end{align}
Hence, using the fact that $1\leq 1+ \sqrt{1-  4\Delta} \leq 2$
we obtain $
\frac{(\sigma^2_1-\sigma^2_2)}{\lambda_1} \le \lambda_{M+N-1}(\mathbf{H}) \le \frac{2(\sigma^2_1-\sigma^2_2)}{\lambda_1}$. 
To obtain the bounds for the condition numbers, by \eqref{LowestTwo_eigH}, we have  
\begin{align*}
 \frac{(M+N)^2}{2 M^2N^2 (\sigma_1^2-\sigma_2^2) } =\frac{\lambda_1^2}{2(\sigma^2_1-\sigma^2_2)}\leq  \kappa(\mathbf{H}) 
 \leq \frac{\lambda_1^2}{(\sigma^2_1-\sigma^2_2)} = \frac{(M+N)^2}{ M^2N^2 (\sigma_1^2-\sigma_2^2) } . 
\end{align*}
{\color{black}{Note the spectral gap of $\mm\Pi^\top \mm\Pi$ is  $\rho(\mm\Pi^\top \mm\Pi)=\sigma_1^2-\sigma_2^2$}}, then this gives the bounds in \eqref{eq:condiH}. 
\end{proof}

\bigskip

The case $M=N$ is of particular interest, and we list the results as a corollary, which follows directly from Theorem \ref{thm:condi_H1}.
\begin{corollary}
\label{cor:condH_M=N}
Let $M=N$ and $\mm\Pi$ be a positive coupling matrix with uniform marginals. The eigenvalues of $\mathbf{H}$ are 
\begin{equation}\label{H-spec_M=N}
\lambda_j(\mathbf{H})= \frac{1}{N}+\sigma_j, \quad
\lambda_{2N+1-j}(\mathbf{H})= \frac{1}{N}-\sigma_{j}, \quad 1\leq j\leq N,  
\end{equation}
where $\{\sigma_j\}$ are the singular values of $\mm\Pi$ in descending order. In particular, $\sigma_1= \frac{1}{N}$, $\lambda_1(\mathbf{H})= \frac{2}{N}$ and $\lambda_{2N}=0$.  
The condition number of $H(\mm\Pi)$ is bounded 
by 
\begin{equation}\label{eq:condiH_M=N}
\frac{2}{N^2 \rho(\mm\Pi^\top \mm\Pi) }  \leq  \kappa(\mathbf{H}) 
 \leq \frac{4}{N^2 \rho(\mm\Pi^\top \mm\Pi) }, 
\end{equation}
{\color{black}{where $\rho(\mm\Pi^\top \mm\Pi)$ is spectral gap of the matrix $\mm\Pi^\top \mm\Pi$.}} 
\end{corollary}

The Sinkhorn algorithm alternatively re-scales the rows and columns of the coupling matrix to achieve the marginal constraints. 
{\color{black}{It produces}} a sequence of coupling matrices $\{\mm\Pi^{(l)}\}$ that converges to $\mm\Pi^*$ entry-wisely, i.e., $\lim_{l\rightarrow +\infty} \Pi_{ij}^{(l)}=\Pi_{ij}^*$ \cite{sinkhorn1964relationship, sinkhorn1967concerning,marshall1968scaling}.  In practice, the Sinkhorn iteration stops when a criterion is met. One stopping criterion is that the marginal distributions of $\mm\Pi^{(l)}$ are entry-wise $\delta$ away from the given $\mm \mu$ and $\mm \nu$. 
Thus, an important question is whether the condition number of $\mathbf{H}$ of $\mm\Pi^{(l)}$ is controlled.

The next proposition shows that if  $\delta = \|\mm\Pi-\mm\Pi^*\|_{F}$ is small with $\|\cdot\|_{F}$ denotes the Frobenius norm, the condition number of $\mathbf{H}$ of $\mm\Pi^{(l)}$ is almost as large as the condition number of $\mathbf{H}^*$. 
The proof is based on Weyl's inequality, and we postpone it to supplementary materials. 

\begin{proposition}[\textbf{Condition number of H matrices in Sinkhorn}] 
\label{prop:condi_Sinkhorn}
Let $\mm\Pi^*$, with uniform marginal distributions, be the optimal coupling matrix minimizing the EOT distance. 
Assume that the coupling matrix $\widehat{\mm\Pi}$ is computed by an early-stopped Sinkhorn algorithm that satisfies 
\begin{equation}\label{eq:assum_perturb}
 \max_{1\le i\le M}\left|\sum_{j=1}^N \widehat{\Pi}_{ij} -\frac{1}{M} \right |\le \delta, \quad \max_{1\le j\le N}\left|\sum_{i=1}^M \widehat{\Pi}_{ij}-\frac{1}{N}\right|\le \delta, \quad \sum_{i,j}  | \widehat{\Pi}_{ij}- \Pi_{ij}^*|^2\leq \delta_2^2. 
\end{equation}
Then, the eigenvalues of 
 $\mathbf{H}$ of $\widehat{\mm\Pi}$ satisfies 
\begin{equation}\label{eq:eig_perturb}
|\lambda_k(\mathbf{H})- \lambda_k(\mathbf{H}^*) |\leq  \delta+\delta_2, \quad 1\leq k\leq N+M. 
\end{equation}
In particular, if $\delta+\delta_2 =t \frac{MN}{M+N}\rho(\mm\Pi^\top \mm\Pi)$ with $t\in[0,1)$, where $\sigma_1, \sigma_2$ denoting the largest two singular values of $\mm\Pi^*$, 
 the condition number of $\mathbf{H}$ of $\widehat{\mm\Pi}$ is bounded by
\[
\frac{1-t\Delta}{(2+t)\Delta} \leq \kappa(\mathbf{H}) 
\leq \frac{1+t\Delta}{(1-t)\Delta}, 
\]
where $\Delta = (\frac{MN}{M+N})^2 \rho(\mm\Pi^\top \mm\Pi)$, while $\frac{1}{2\Delta} \leq \kappa(\mathbf{H}^*) \leq \frac{1}{\Delta}$. 
\end{proposition}

The above bounds apply to general $\mathbf{H}$ matrices in entropy-regularized Sinkhorn algorithms. However, these bounds do not show explicit dependence on $\epsilon$, the strength of regularization. In the next section, we study the dependence of the condition number on $\epsilon$ and $N$ for point-clouds datasets.

\subsection{Ill-conditioned H matrices from data clouds}\label{subsec:NumericalSpecStudy}
We investigate in this section the condition number of the $\mathbf{H}^*$ for a specific example of EOT that matches data clouds with $N$ points. In this simple setting, $M=N$ and both marginals are uniform, so the condition number of the $\mathbf{H}^*$ matrix is 
$\kappa =\frac{2}{N\lambda_{2N-1}}$ {\color{black}{by Corollary}}~\ref{cor:condH_M=N}. 
Therefore, it suffices to investigate the dependence of the smallest positive eigenvalue, $\lambda_{2N-1}$, on $N$ and $\epsilon$.

We show that the smallest positive eigenvalue of the $\mathbf{H}^*$ matrix can decay at rate $O(e^{-\frac{1}{\epsilon}})$ for a fixed $N$ and at $O(1/N)$ for a fixed $\epsilon$. These asymptotic results are proved for equally-spaced points on the unit circle in Example \ref{exp:H-ill-conditioned}, and are numerically demonstrated for random data clouds sampled from a uniform distribution. {\color{black}{So the $\mathbf{H}^*$ matrix in practice could be very ill-conditioned, and the proper regularization to resolve numerical instability is necessary. }}

\begin{example}[Equally spaced points on the unit circle] 
\label{exp:H-ill-conditioned}
Consider $N$ equally spaced points on the unit circle $\{ \mathbf{y}_i= \bcm \cos x_i & \sin x_i \ecm\}_{i=0}^{N-1}$, where $x_i =\frac{2\pi i}{N}$. Let $\mm\mu=\frac{1}{N}\mathbb{1}_N$ be the uniform distribution, we are interested in the spectrum of $\mathbf{H}^*$ associated to the symmetric entropic regularized optimal transport loss $\text{OT}_\epsilon(\mathbf{C},  \mm\mu, \mm \mu)$. 
The coupling matrix is 
$$
\mm\Pi^*= \argmin{\mm\Pi\in \mathbf{U}(\mm \mu, \mm \mu) }  \sum_{i=1}^N\sum_{j=1}^N C_{ij} \Pi_{ij} + \epsilon \text{KL}(\mm\Pi, {\mm \mu} \otimes \mm\mu) 
$$ 
with $\epsilon>0$, where $C_{ij}=\|\mathbf{y}_i - \mathbf{y}_{j}\|_2^2$. Denote its condition number by $\kappa(\mathbf{H}^*)= \frac{\lambda_1(\mathbf{H}^*)}{\lambda_{2N-1}(\mathbf{H}^*)}$. 
Then, the following statements hold true. 
\begin{itemize}
   \item[(a)] $\mm\Pi^*=\frac{\mathbf{K}}{\lambda_1(\mathbf{K}) N} $, where the Gibbs kernel $\mathbf{K}\in \R^{N\times N}$ is a symmetric matrix with entries $K_{ij}=\exp\left(-\frac{C_{ij}}{\epsilon}\right)$ and $\lambda_1(\mathbf{K})$ is the largest eigenvalue of $\mathbf{K}$. 
   \item[(b)] The first two singular values of $\mm\Pi^*$ are $\sigma_1(\mm\Pi^*) = \frac{1}{N}$, $ \sigma_2(\mm\Pi^*) = \frac{\lambda_2(\mathbf{K})}{\lambda_1(\mathbf{K})N}$, and the largest and smallest positive eigenvalues of $\mathbf{H}^*$ are $\lambda_1(\mathbf{H}^*) = \frac{2}{N}$, and $\lambda_{2N-1}(\mathbf{H}^*) = \frac{1}{N}- \sigma_2(\mm\Pi^*)$, where $\lambda_2(\mathbf{K})$ is the second largest eigenvalue of $\mathbf{K}$. 
   
   \item[(c)] The smallest positive eigenvalue of $\mathbf{H}^*$ and the condition number of $\mathbf{H}^*$ satisfies 
		\begin{align}
			\lim_{\epsilon\rightarrow 0^+}\lim_{N\rightarrow +\infty}\frac{N\cdot \lambda_{2N-1}(\mathbf{H}^*) }{\epsilon}  &=\frac{1}{4},  \quad  & \lim_{\epsilon\rightarrow 0^+}\lim_{N\rightarrow +\infty}\epsilon\cdot\kappa(\mathbf{H}^*)  & = 8,  \label{eq:limit_N_then_e} \\
					\lim_{N\rightarrow+\infty}\lim_{\epsilon\rightarrow 0^+} \frac{ N \lambda_{2N-1}(\mathbf{H}^*) } {r_{N,\epsilon}} & = 4\pi^2,   \quad &
				\lim_{N\rightarrow +\infty}\lim_{\epsilon\rightarrow 0^+ } r_{N,\epsilon} \kappa(\mathbf{H}^*)  & = \frac{1}{2\pi^2}, \label{eq:limit_e_then_N}
			\end{align}
			where  $r_{N,\epsilon} = N^{-2}\exp\left(-\frac{4\sin^2(\pi/N)}{\epsilon }\right)$. 
\end{itemize}
We postpone the derivation to supplementary materials.
\end{example}

\begin{figure}[tp!]
\centering
\vspace{-0.5 cm}\hspace{-1mm}
\subfigure[]{\includegraphics[width =0.48 \textwidth]{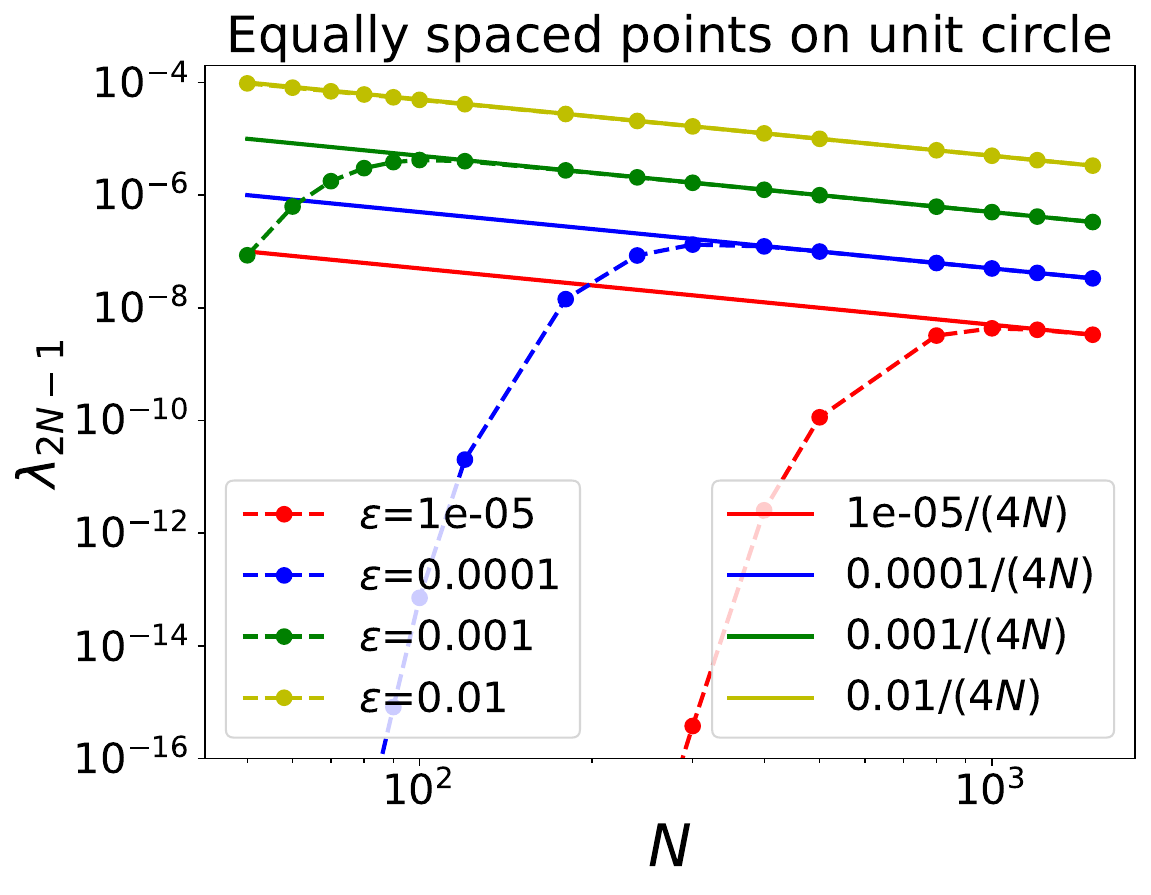}}
\hspace{-2mm}
\subfigure[]{\includegraphics[width =0.48 \textwidth]
{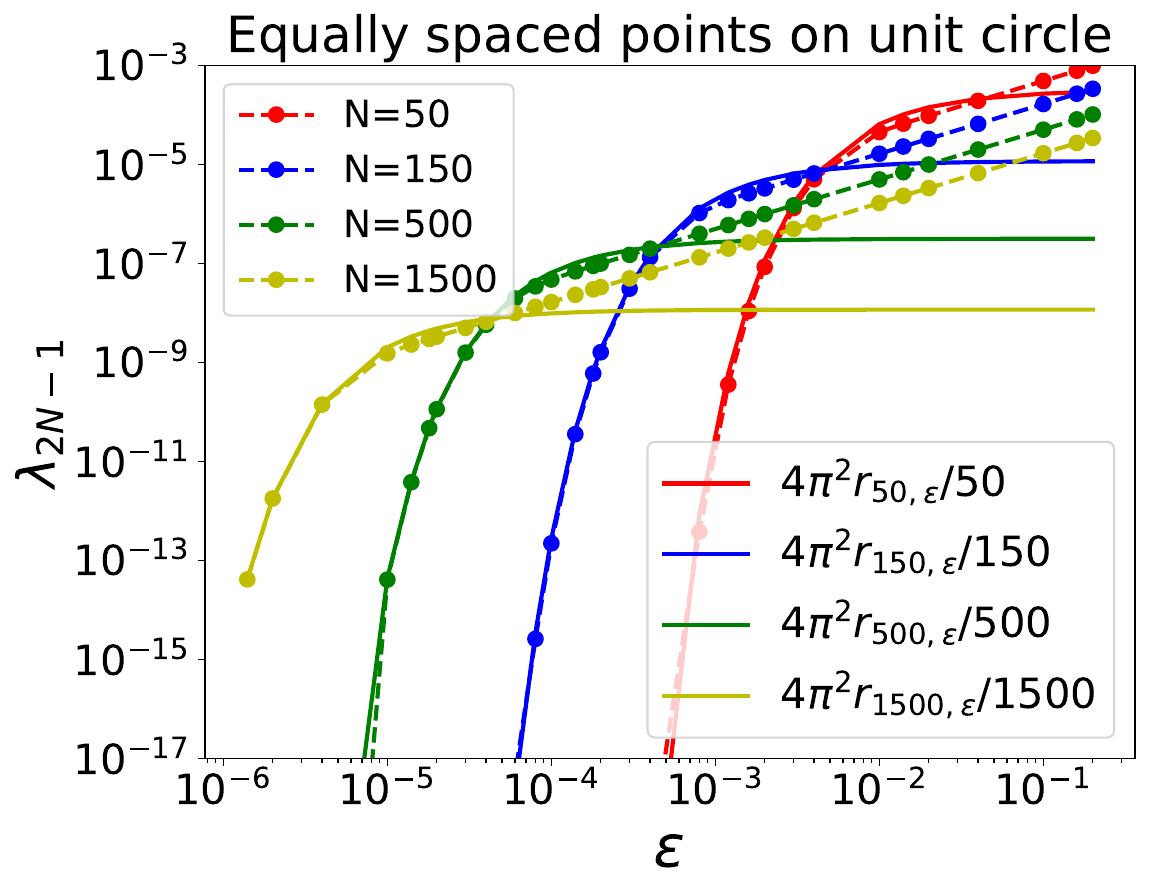}}

\hspace{-2mm}
\subfigure[]{\includegraphics[width =0.48\textwidth]{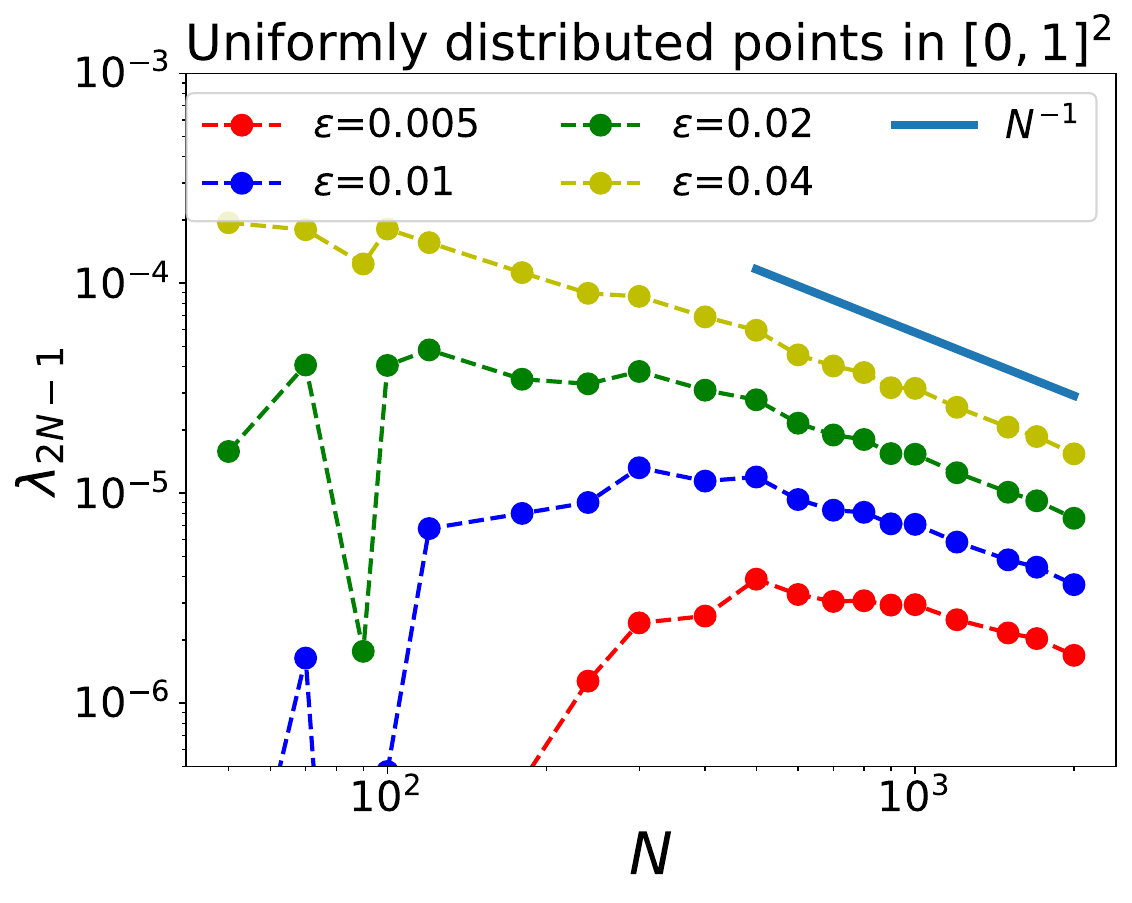}}
\hspace{-1mm}
\subfigure[]{\includegraphics[width =0.48 \textwidth]{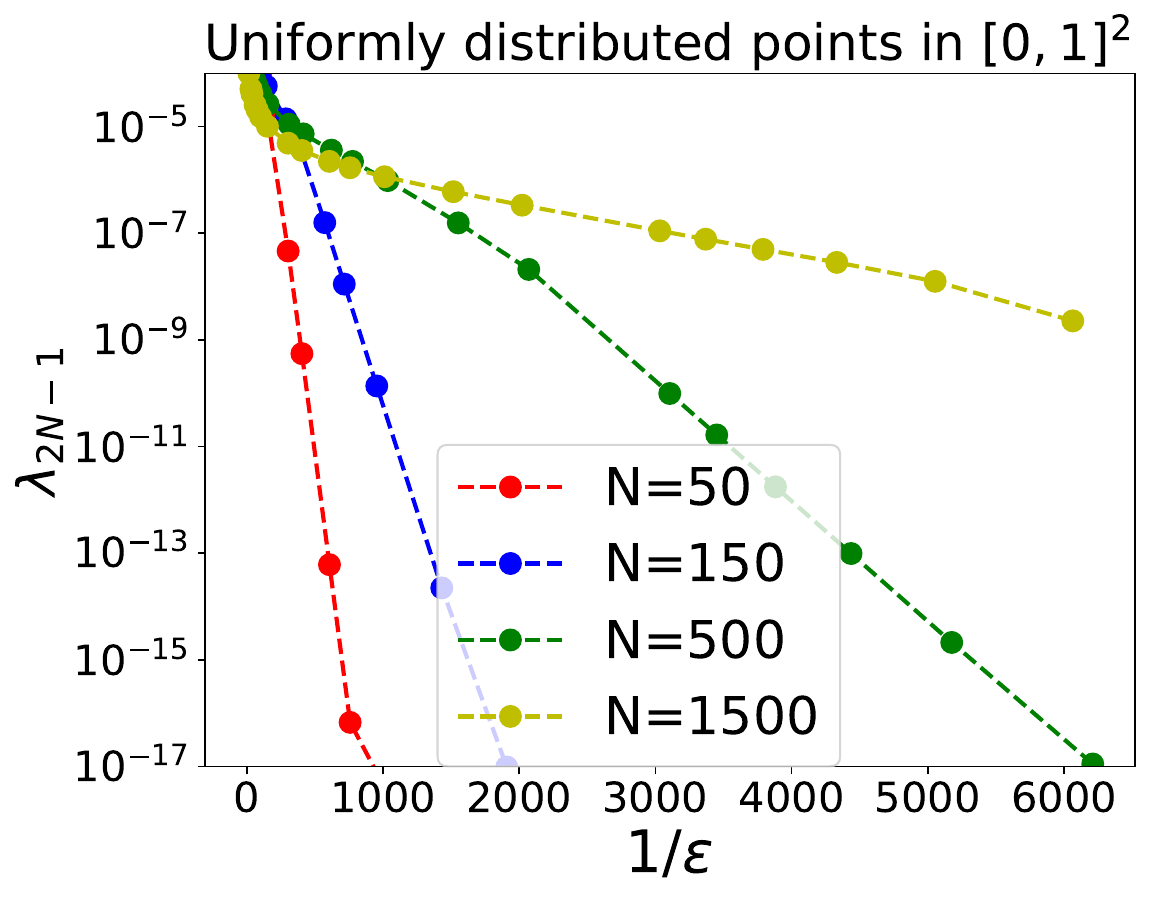}}
\vspace{-0.4 cm}
\caption{ Decay of the smallest positive eigenvalue $\lambda_{2N-1}$ in $N$ and $\epsilon$.  Equally spaced points on the unique circle: \textbf{(a)} $\lambda_{2N-1}\approx  \frac{\epsilon}{4N}$ when 
$N > \frac{2\pi}{\sqrt{\epsilon}}$; 
\textbf{(b)}  $\lambda_{2N-1}\approx 4\pi^2 r_{N,\epsilon} $ when $\epsilon< \frac{4\pi^2}{N^2}$. Uniformly distributed point cloud in unit square $[0,1]^2$: \textbf{(c)} $\lambda_{2N-1}=O(\frac{1}{N})$ when $N$ is large; \textbf{(d)} $\lambda_{2N-1}=O(e^{-\frac{1}{\epsilon}})$ when $\epsilon$ is small.
 }
\label{Fig: Data_Eps_spect}
\vspace{-0.7cm}
\end{figure}

  For a fixed $\epsilon$, when $N$ is large enough, the smallest positive eigenvalue $\lambda_{2N-1}$ scales as $\frac{\epsilon}{4N}$, which is numerically illustrated in Figure \ref{Fig: Data_Eps_spect}(a).
  Then the condition number scales as $\frac{8}{\epsilon}$. 
  Meanwhile for fixed $N$, when $\epsilon$ is small enough (e.g., when $\epsilon<\frac{4\pi^2}{N^2}$), the smallest positive eigenvalue scales as $\frac{4\pi^2 r_{N,\epsilon}}{N}$, which is numerically illustrated in Figure \ref{Fig: Data_Eps_spect}(b). 
  Then the condition number scales as $\frac{1}{2\pi^2 r_{N,\epsilon}}$, which grows  exponentially at rate $O(e^{-\frac{1}{\epsilon}})$. 
  For example, for $N=50$ and $\epsilon=0.0001$,  the condition number of $\mathbf{H^*}$ is larger than $10^{70}$; in this case, a truncated SVD for $\mathbf{H}^*$ is crucial in calculating the Hessian of EOT and the gradient of Sinkhorn distance. On the other hand, when $N=1500$ and $\epsilon=0.0001$, the condition number of $\mathbf{H}^*$ is only about $8\times 10^4$.
In addition, this scale phenomenon is also observed for some random datasets as well.

Next, we further numerically investigate the case of point-clouds datasets that are sampled from  the uniform distribution in the unit square $[0,1]^2$.  Similarly, we are interested in the spectrum of $\mathbf{H}^*$ associated with the symmetric $\text{OT}_\epsilon(\mathbf{C},  \mm\mu, \mm \mu)$.
In Figure \ref{Fig: Data_Eps_spect}(d) show that $\lambda_{2N-1} = O(e^{-\frac{1}{\epsilon}})$ when $\epsilon$ is small enough for each fixed N, and Figure \ref{Fig: Data_Eps_spect}(c) show $\lambda_{2N-1} = O(\frac{1}{N})$ when $N$ is large for each fixed $\epsilon$. These asymptotic orders are the same as the analytical results proved in Example \ref{exp:H-ill-conditioned}, but the exact limit depends on the distribution of data points, and it is beyond the scope of this study.  

In summary, $\mathbf{H}^*$ can be severely ill-conditioned with the smallest eigenvalue at the order of $O(e^{-\frac{1}{\epsilon}})$ when $\epsilon$ is small, or $O(\frac{1}{N})$ when $N$ is large. Thus, when solving a linear system with $\mathbf{H}^*$, it is important to properly regularize the ill-posed inverse problem.

\section{Hessian computation: runtime, accuracy, and success rate} \label{sec:Hess_comp}
Our analytical approach enables efficient and accurate computation of the Hessian matrix. Here we compare it with the current two state-of-the-art approaches suggested by \textit{OTT}: \textit{unroll} and \textit{implicit differentiation}. The details on these approaches are discussed in Section \ref{sec:prev}.

We use the point-cloud datasets sampled from the uniform distribution in unit square again. 
The task is to calculate the Hessian tensor $\mathcal{T}$ of 
$\text{OT}_\epsilon(\mathbf{C},  \mm\mu, \mm\mu)$ respect to the source data $\mathbf{Y}$, where $\mm\mu=\frac{1}{N}\mathbb{1}_N$. By proposition \ref{proposition:Hessian_sum}, the Hessian satisfies the marginal identity, i.e., $\sum_{k=1}^M    \mathcal{T}_{k\cdot s\cdot } = 2\mu_s \mathbb{I}_d$. We evaluate the accuracy of the computed Hessian by the  marginal error: 
\begin{align}\label{MargError}
    \text{error} = \sum_{t\ne l}(\sum_k\mathcal{T}_{ktsl})^2  + \sum_{t} (\sum_k\mathcal{T}_{ktst}-2\mu_s)^2. 
\end{align}
All simulations are performed on a single Nvidia A100 GPU using double-precision. The threshold $\alpha$ in truncated SVD is set to $\alpha=10^{-10}$. 

\begin{figure}[htp!]
    \centering
{\includegraphics[width =0.8 \textwidth]{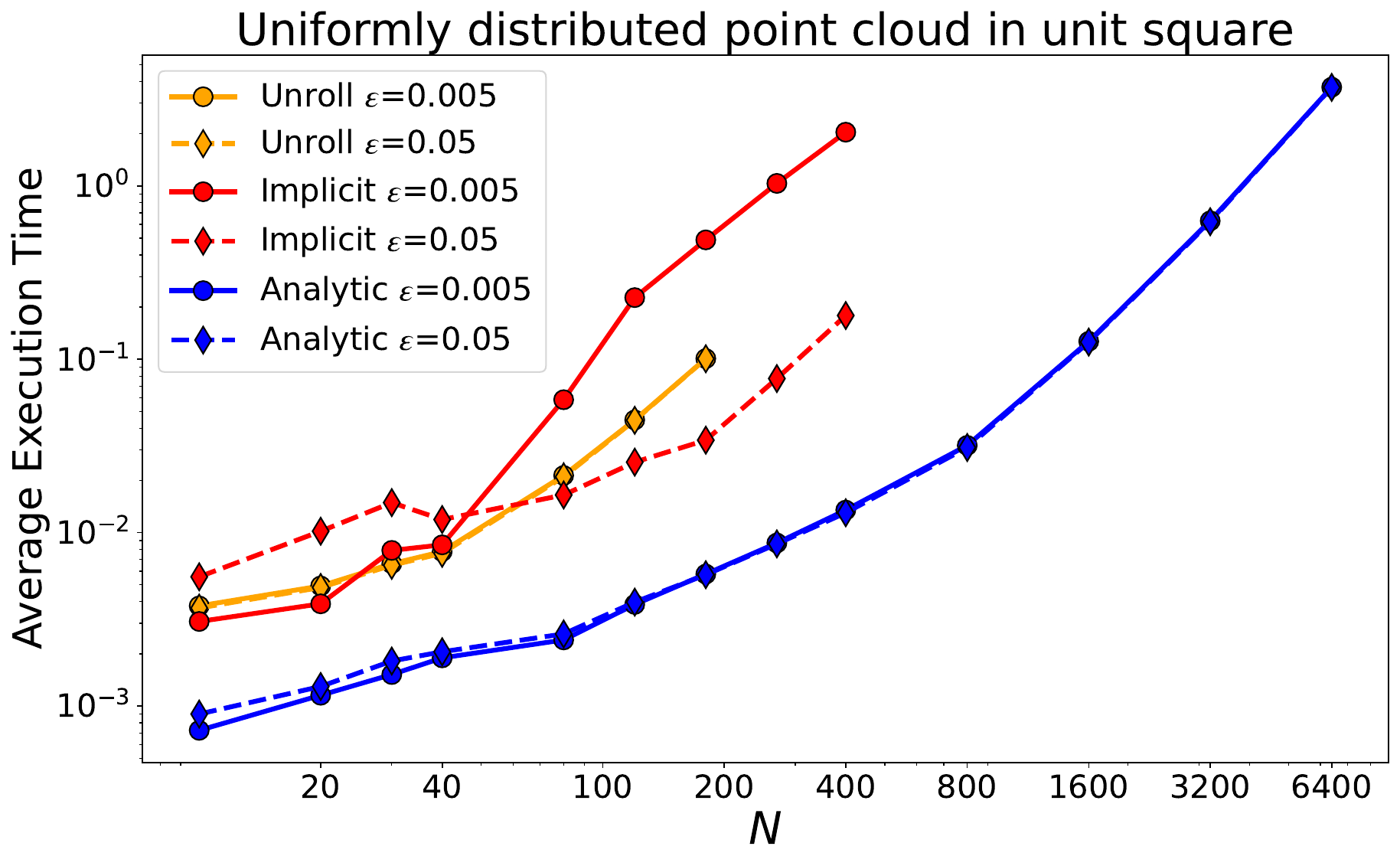}}

{\includegraphics[width =0.84 \textwidth]{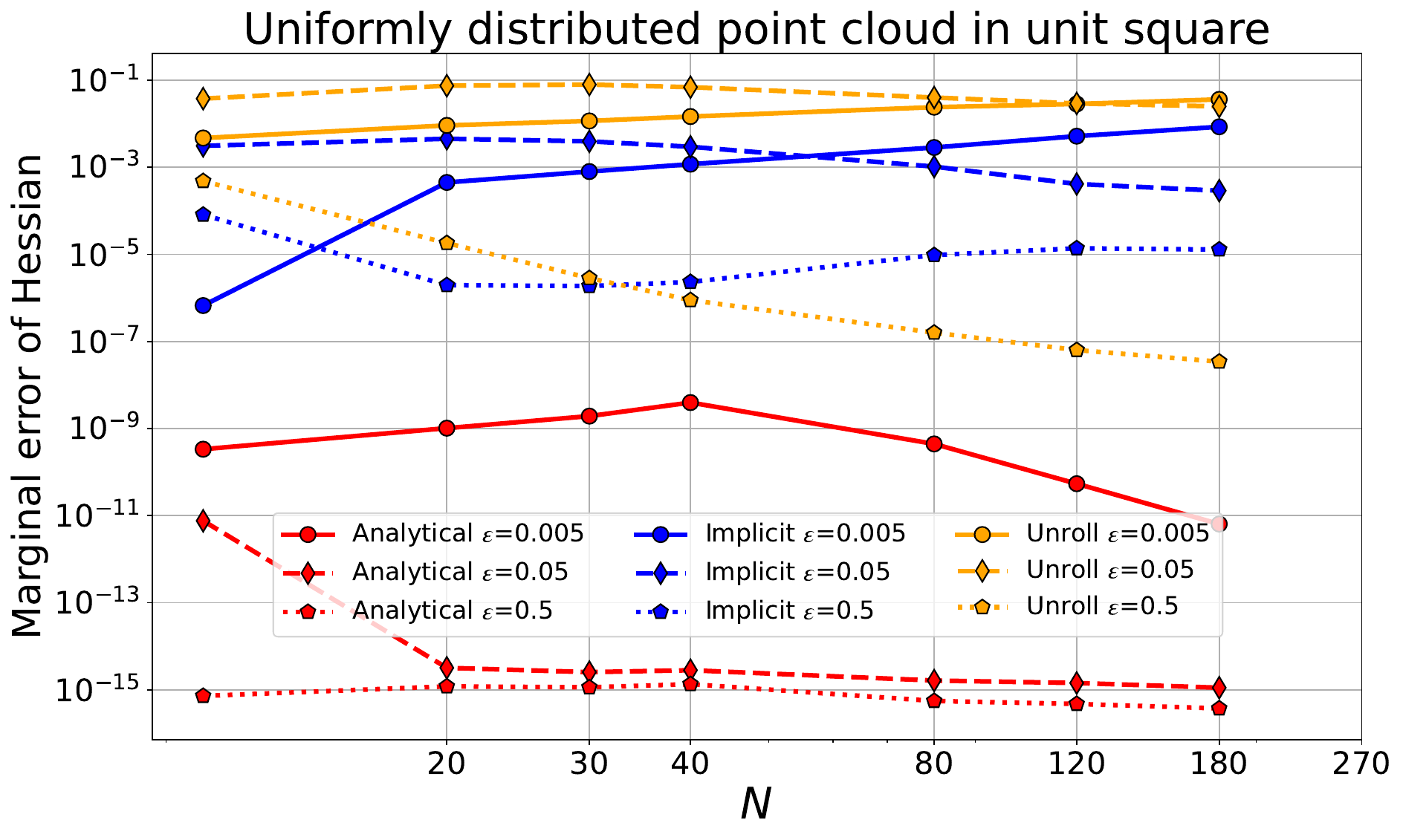}}
    \caption{Comparison of runtime (in seconds) and marginal error for  Hessian computing $\frac{d^2\text{OT}_\epsilon(\mathbf{C}, \mm \mu,  \mm \mu)}{d\mathbf{Y}^2}$ among three approaches: unroll, implicit differentiation and analytic expression with regularization (ours). }
\label{fig:Hessian_accuracy_speed}
\vspace{-0.7cm}
\end{figure}

\textbf{Runtime.} Figure \ref{fig:Hessian_accuracy_speed}(a) shows the average execution time in 10 independent tests for the three approaches with $N\in [10, 6400]$ and $\epsilon\in \{0.005, 0.05\}$, corresponding to low and median regularization regimes. The unrolling and implicit differentiation approaches fail in all 10 tests due to insufficient memory when $N>180$ and $N>400$, respectively. However, our analytical approach remains effective for all $N$, even beyond $N=5000$. Additionally, when all three approaches work, our analytical approach is faster by at least one order of magnitude.

\textbf{Accuracy.}  Figure \ref{fig:Hessian_accuracy_speed}(b) shows the 
average marginal error of the Hessian computed by the three approaches in 100 independent tests. Here we consider $\epsilon\in \{0.005, 0.05\}$ and $N\in[10,180]$ where all three approaches work.  
Both implicit differentiation and unrolling approaches perform poorly across all parameter settings. In contrast, our analytical approach is significantly more accurate by 3-8 orders of magnitude.

\textbf{Success rate.} Table \ref{tab:methods_comparison} further highlights the reliability of our analytical approach and the importance of regularization by reporting the success rate in  100 independent tests. A test is considered successful if the marginal error of the Hessian \eqref{MargError} is less than 0.1.  In the most singular parameter setting, $N=10$ and $\epsilon=0.005$ (as discussed in section \ref{subsec:NumericalSpecStudy}),  the implicit differentiation approach fails 97\% of the tests due to numerical instability.  Importantly, if we do not regularize the problem using truncated SVD and instead apply the least square solver directly to solve the linear system, the analytical approach results in large errors ranging from $10^{-7}$ to $10^{2}$ in 15\% of the tests. Therefore, proper regularization is crucial when the problem is ill-posed.

\begin{table}[ht]
\centering
\begin{tabular}{lcccc}
\toprule
            Method &  Unroll &  Implicit & Analytical(no reg) &    Analytical(with reg)\\
\midrule
            $N=10$ &  0.78 &  0.03 &   0.85 &     1.00 \\
          $N=20$ &  0.68 &  0.18 &   0.99 &     1.00 \\
$N=120$ &  0.00 &  1.00 &   1.00 &     1.00 \\
  $N=1600$ &  0.00 &  0.00 &   1.00 &     1.00 \\
\bottomrule
\end{tabular}
\caption{Success rates of the three approaches for $N\in\{10,20,120,1600\}$ and $\epsilon=0.005$. A test is called successful if the marginal error of Hessian \eqref{MargError} is less than 0.1.}
\label{tab:methods_comparison}
\vspace{-0.5cm}
\end{table}
To conclude, our analytical approach with regularization enables efficient and accurate computation of the Hessian of EOT, significantly outperforming other current state-of-the-art approaches by a large margin in terms of runtime, accuracy and success rate.  


\section{Applications to Shuffled Regression}
\label{sec:Application}
In this section, we apply our proposed algorithms to solve shuffled regression problems.
It is formulated as the multivariate regression model $\mathbf{y}^*= \mathbf{x}\theta+ \mm\xi$, where $\mathbf{x}\in \mathbb{R}^D, \mathbf{y}^*\in \mathbb{R}^d, \theta\in \mathbb{R}^{D\times d}$ and $\mm\xi$ is the Gaussian noise independent of $\mathbf{x}$. The correspondence between $(\mathbf{X}, \mathbf{Y}^*)$ is missing. 
Our goal is to estimate the optimal $\theta^*$ using EOT distance as the loss function in the \eqref{framework_prob}. {\color{black}{We choose both probability weights $\mm\mu$ and $\mm \nu$ as uniform weights.}}
This approach generalizes to datasets without requiring $\mathbf{X}$ and $\mathbf{Y}^*$ to have the same number of rows. The gradient and the Hessian of the EOT distance with respect to the parameters $\theta$ are simplified as 
\begin{align}\label{linear_diff_EOT}
   & \frac{d \text{OT}_\epsilon(\mathbf{C}_\theta, \mm \mu, \mm \nu) }{d \theta} =  \mathbf{X}^\top \frac{d\text{OT}_\epsilon(\mathbf{C}_\theta, \mm \mu, \mm \nu)}{d \mathbf{Y}}, \nonumber \\&\left(\frac{d^2\text{OT}_\epsilon(\mathbf{C}_\theta, \mm \mu, \mm \nu) }{d\theta^2}\right)_{mtnl} = \sum_{k=1}^M \sum_{s=1}^M \mathbf{X}_{sm} \mathcal{T}_{ktsl}\mathbf{X}_{kn}
\end{align}
for $t,l = 1, \dots, d$ and $m,n=1,\dots, D$. The EOT distance is generally non-convex with respect to $\theta$, so optimization may not converge to the optimal $\theta^*$. Our focus is on the convergence speed to a local minimum. First-order methods may converge to a local minimum but require many iterations due to the complicated landscape of the loss function.  To accelerate optimization, we propose a two-stage approach. First, we use stochastic gradient descent (SGD) with a random subset of 
$\mathbf{X}$ and the full batch of $\mathbf{Y}^*$ to quickly approach the local minimum. Then, we switch to a relaxed Newton's method, using the updated parameter $\hat{\theta}$ as the initial condition. The relaxed-Newton's method uses step-size $\gamma<1$. In practise, we switch from SGD to relaxed-Newton when the computed Hessian $\frac{d^2\text{OT}_\epsilon(\mathbf{C}_\theta, \mm \mu, \mm \nu)}{d\theta^2}$ is positive definite. 
The algorithm is summarized in Algorithm~\ref{Alg_opt_prob}. {\color{black}{In the situation that the Hessian is too expansive to calculate multiple times, we can use the Hessian calculated from symmetric entropic OT problem of $\mathbf{Y}^*$ to $\mathbf{Y}^*$,  
as the approximation of the true Hessian in Newton's method.} }

\begin{algorithm}
\caption{Two-stage algorithm to estimate optimal $\theta^*$ of EOT distance \eqref{reg_OT_problem}.
\label{Alg_opt_prob}}
\textbf{Input:} Data $\mm X$, target data $\mm{Y}^*$, entropy regularization strength $\epsilon$, truncated SVD threshold $\alpha$; initial guess of $\theta^{(0)}$, SGD learning rate $r_s$, mini batch size $n_s$, maximum epochs $T$; Relaxed Newton learning rate $r_n$. \\
\textbf{Output:} Estimated optimal $\theta^*$, regularized optimal transport loss $\text{OT}_\epsilon(C_{\theta^*},\mm{\mu}, \mm{\nu})$. 
\begin{algorithmic}[1]
\State Set $\mm{\tilde{\mu}}\leftarrow \frac{1}{n_s}\mathbb{1}_{n_s}$, $ \mm \mu \leftarrow \frac{1}{N}\mathbb{1}_{N}$ and $ \mm \nu \leftarrow\frac{1}{N}\mathbb{1}_{N}$.

\noindent {\bf Stage 1 SGD }:

\For{$t=0, \dots, T-1$}, 
\State Randomly sample $n_s$ rows of $\mm X$, denote as $\mm{\tilde{X}}$. 

\State Compute $\mm{\tilde{Y}}\leftarrow \mm{\tilde{X}}\cdot \theta^{(t)}$ and  cost matrix $C^{\mm{\tilde{Y}} \rightarrow \mm Y^*}_{ij}\leftarrow   \|\mm{\tilde{y}}_i - \mm y^*_j\|_2^2$.

\State $\theta^{(t+1)}\leftarrow \theta^{(t)}- r_s\mm{\tilde{X}}^\top \frac{d\text{OT}_\epsilon(C^{\mm{\tilde{Y}} \rightarrow \mm Y^*},\mm{\tilde{\mu}} ,\mm \nu)}{d \mm Y} $

\If{ $\frac{d^2\text{OT}_\epsilon(C_\theta, \mm \mu, \mm \nu) }{d\theta^2}|_{\theta = \theta^{(t+1)}}$ is positive definite }
\State stop Stage 1 with the current $\hat\theta \leftarrow \theta^{(t+1)}$.
\EndIf

\EndFor

\noindent {\bf Stage 2 Relaxed Newton's method }: 

\State Set $\theta^{(0)}\leftarrow \hat\theta$. 

\For{$t=0, \dots, T-1$}
\State Compute $\mm{Y}\leftarrow \mm{X}\cdot \theta^{(t)}$ and  cost matrix $C_{ij}\leftarrow \|\mm{y}_i - \mm y^*_j\|_2^2$

\State 
$\theta^{(t+1)}\leftarrow \theta^{(t)}- r_n \left(\frac{d^2\text{OT}_\epsilon(C_{\theta^{(t)}}, \mm \mu, \mm \nu) }{d\theta^2} \right)^{-1} \left(\mm{\tilde{X}}^\top \frac{d\text{OT}_\epsilon(C_{\theta^{(t)}},\mm{\mu} ,\mm \nu)}{d \mm Y}\right) $
\If{$\text{OT}_\epsilon(C_{\theta^{(t+1)}},\mm{\mu} ,\mm \nu)$ doesn't improve}
\State Quit Stage 2 with $\theta^*\leftarrow \theta^{(t+1)}$.
\EndIf

\EndFor

\end{algorithmic}
\end{algorithm}

\subsection{Shuffled Regression with Gaussian Mixtures}\label{subsec:ShuffReg}

We first generate $N=500$ data points $\mathbf{X}\in \mathbb{R}^5$ from a Gaussian mixture distribution with three clusters. {\color{black}{The mean vectors for these clusters are randomly chosen, and the covariance matrices are 
 $0.3^2\mathbb{I}_5, 0.05^2\mathbb{I}_5$ and $0.6^2\mathbb{I}_5$, respectively}}.  The parameter $\theta^*\in \mathbb{R}^{5\times 2}$ is generated with components $\theta^*_{mt}\sim \mathcal{N}(0,1)$, and the Gaussian noise $\mm \xi\in \R^2$ follows $\mm\xi \sim \mathcal{N}(\mathbb{0}, 0.04\mathbb{I}_2)$.
We then compute $\mathbf{y}_i^* = \mathbf{x}_i\theta^*+\mm{\xi}_i$, randomly and completely permute the order of $\mathbf{y}_i^*$, removing $\mathbf{X}$-to-$\mathbf{Y}^*$ correspondence. 

\begin{figure}[t]
    \centering
\subfigure[]{\includegraphics[width =0.4 \textwidth]{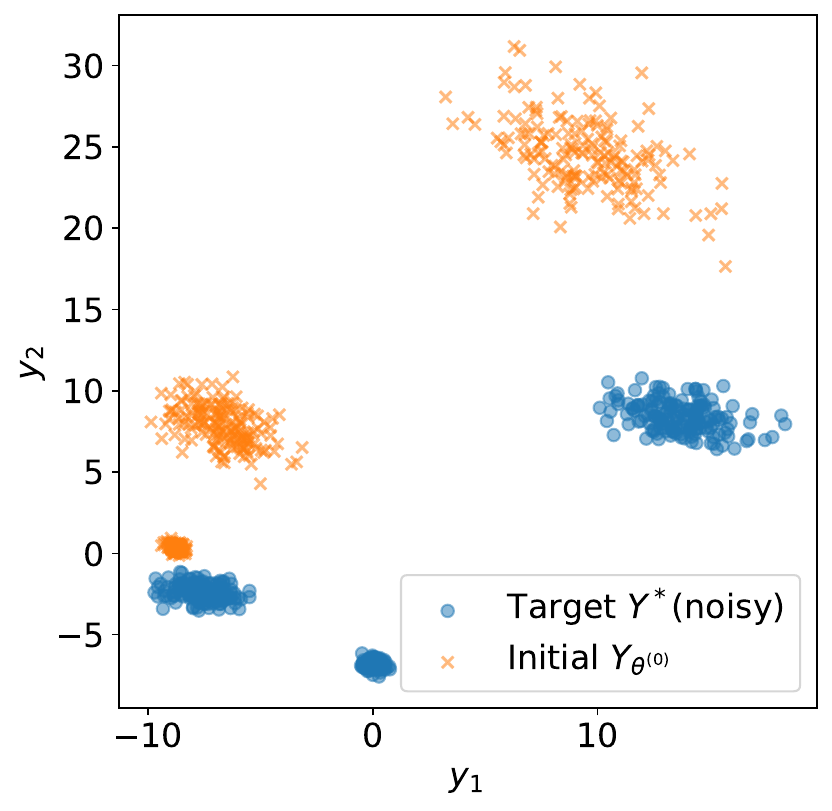}} 
\subfigure[]{\includegraphics[width =0.4 \textwidth]{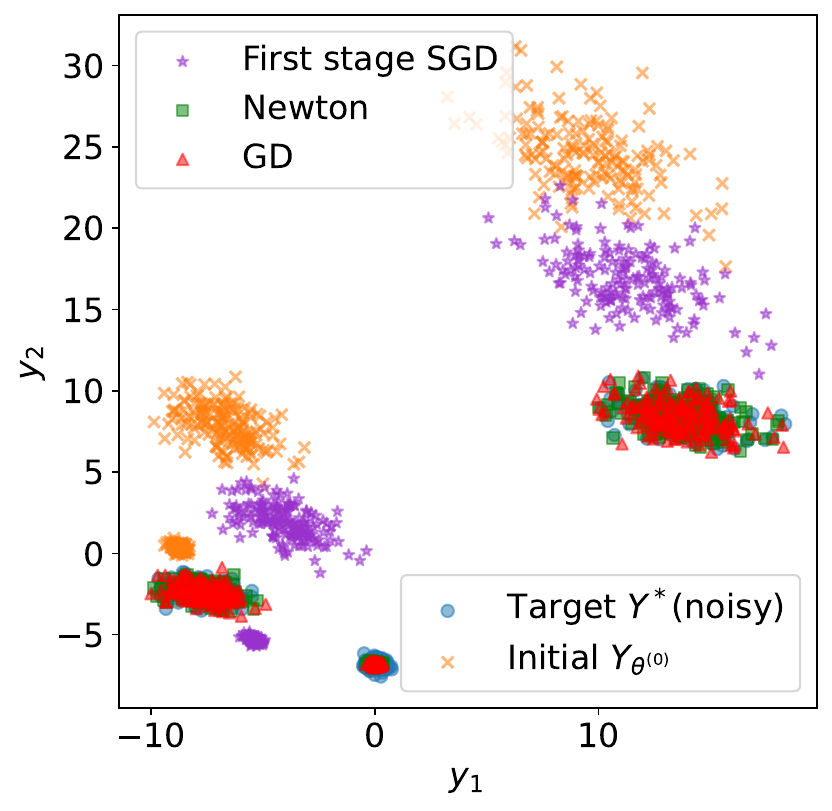}}

\subfigure[]{\includegraphics[width =0.325 \textwidth,valign=t]{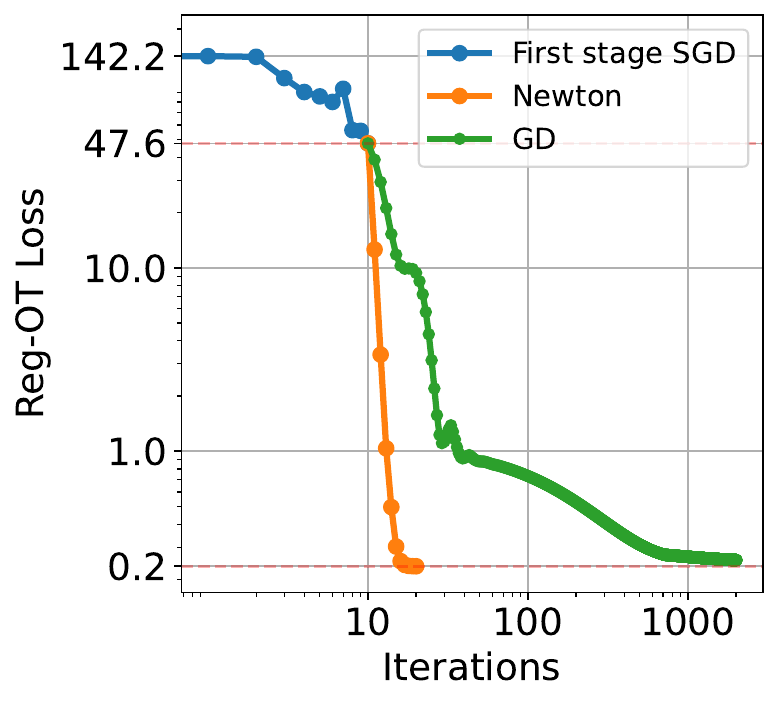}}
\subfigure[]{\includegraphics[width =0.325\textwidth,valign=t]{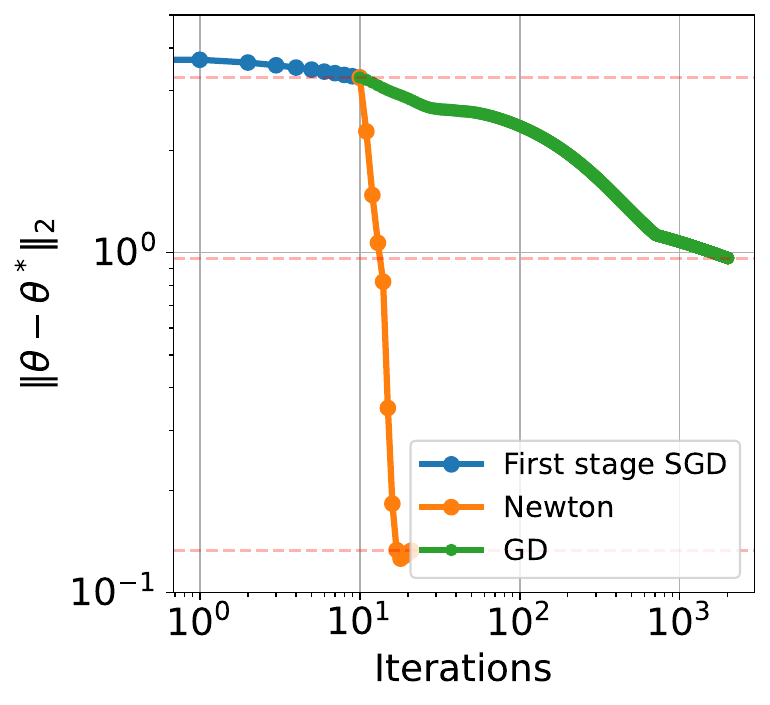}}
\subfigure[]{\includegraphics[width =0.268\textwidth,valign=t]{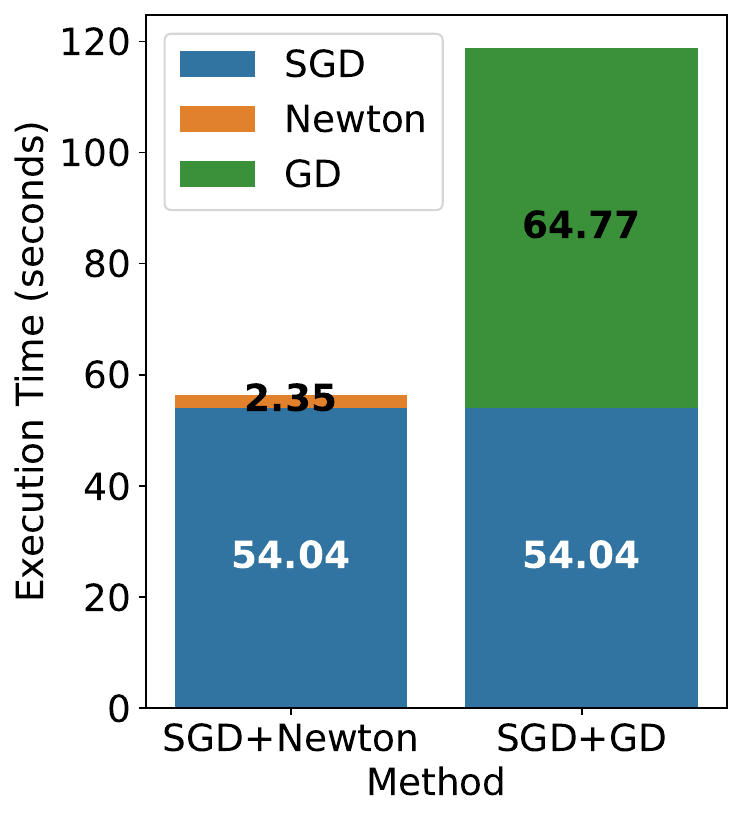}}
    \caption{Shuffled Regression with Gaussian Mixtures.   
    }
    \label{fig:shuffuled_regression}
    \vspace{-0.7cm}
\end{figure}

Starting with an random initial condition $\theta^{(0)}$ from the standard normal distribution, the target data $\mathbf{Y}^*$ and the initial data $\mathbf{Y}_{\theta^{(0)}}$ are shown in Figure \ref{fig:shuffuled_regression}(a).  We use the two-stage algorithm described above. In the first stage, we perform 10 iterations of SGD on 100 random source data points with a learning rate of 0.001. In the second stage, we use a relaxed Newton's method with a learning rate of 0.5. We compare this to a gradient descent (GD) method with a learning rate of 0.001.

Figure \ref{fig:shuffuled_regression}(b) shows that both methods correctly map the data $\mathbf{X}$ to the target data $\mathbf{Y}^*$.
Figure \ref{fig:shuffuled_regression}(c-d) shows that both methods converge to the optimal $\theta^*$,  but the relaxed Newton's method is faster and more accurate.  The relaxed Newton's method converges in 12 iterations with a runtime of 2.35 seconds, while GD takes 2000 iterations and 64.77 seconds, which is 27 times longer. Additionally, the relaxed Newton's method achieves nearly one order of magnitude better accuracy in terms of the $L_2$ error $\|\theta-\theta^*\|_{2}$.  

Further analysis shows that the eigenvalues of Hessian with respect to $\theta^*$ range from $10^{-2}$ to $10^2$, indicating that the optimal parameter lies in a long, narrow, flat valley, causing the gradient descent method to converge slowly.

\subsection{3D Point Cloud Registration}\label{subsec:3Dpoints}
In this section, we extend our method to 3D point clouds registration, a critical task in computer vision. The goal is to find a spatial transformation that aligns two 3D data clouds without knowing the correspondence, known as simultaneous pose and correspondence registration \cite{qi2017pointnet,katageri2024metric}.

Using the MobilNet10 dataset \cite{qi2017pointnet}, we create a study room with a chair (500 points), a desk (1500 points), and a sofa (1500 points), denoted as $\mathbf{X}$. We apply a linear transformation including random rotation and scaling, and add Gaussian noise: $\mathbf{Y}^* = \mathbf{X}\theta^* + \mm \xi$ with $\mm \xi\sim \mathcal{N}(0, 4\times 10^{-4}\mathbb{I}_3)$. The rows of $\mathbf{Y}^*$ are randomly permuted to remove correspondence. 

We use the two-stage algorithm described above. The initial parameter $\theta^{(0)}$
  is a standard Gaussian perturbation of the optimal parameter $\theta^*$. In the first stage, we perform 5 iterations of SGD on 500 random data points with a learning rate of 0.1. In the second stage, we use a relaxed Newton's method with a learning rate of 0.5. For comparison, the GD method uses a learning rate of 0.1.
  
\begin{figure}
\centering

\includegraphics[width=0.75\textwidth]{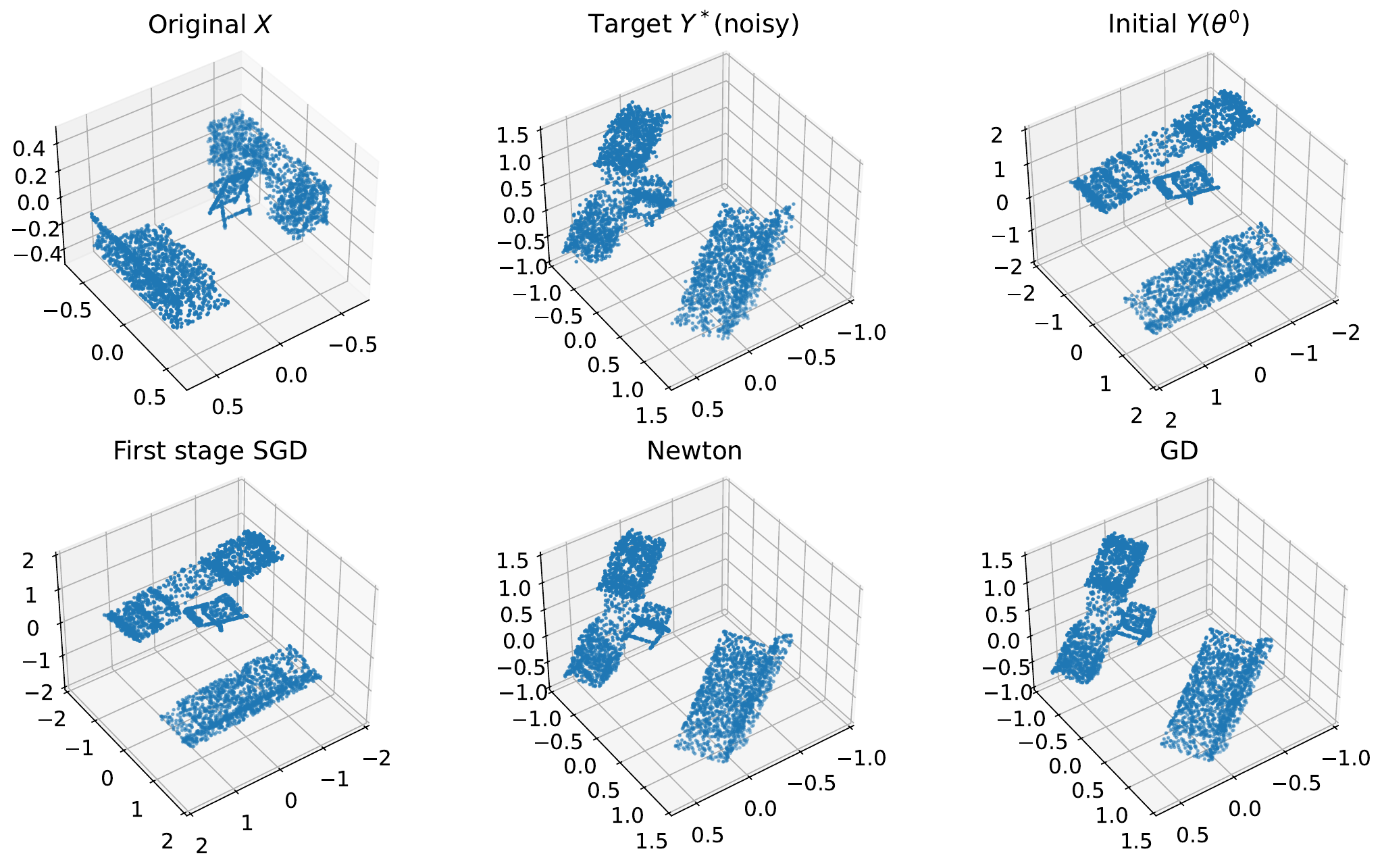}

\includegraphics[width=0.28\textwidth,valign=t]{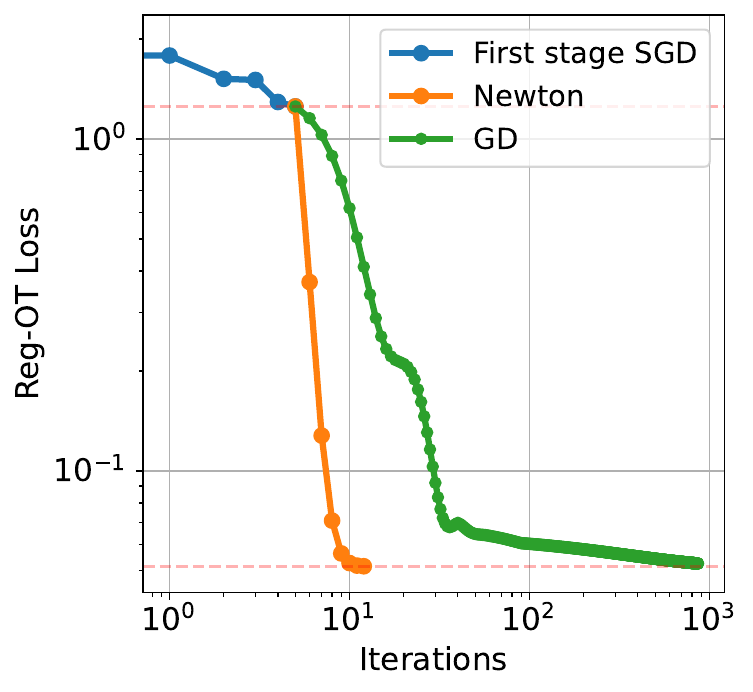}
\includegraphics[width=0.28\textwidth,valign=t]{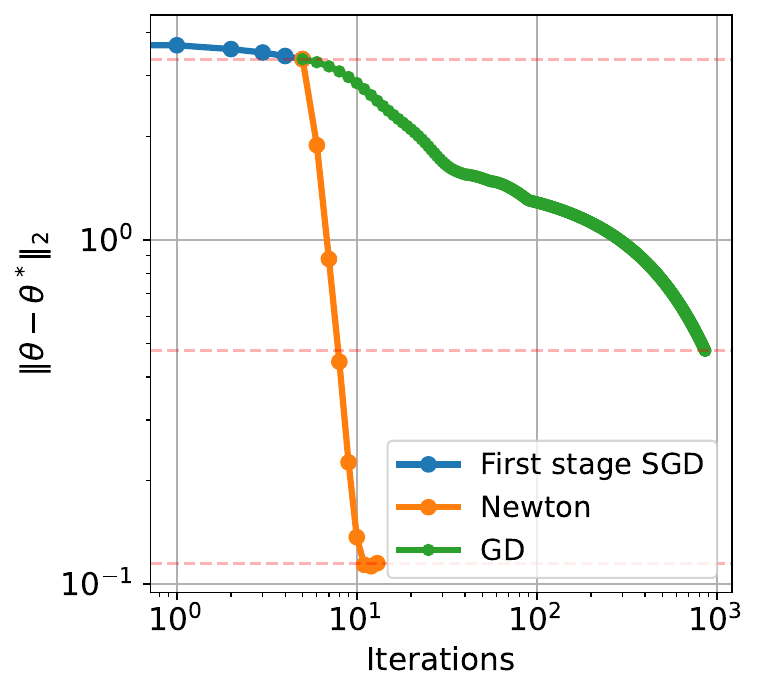}
\includegraphics[width=0.226\textwidth,valign=t]{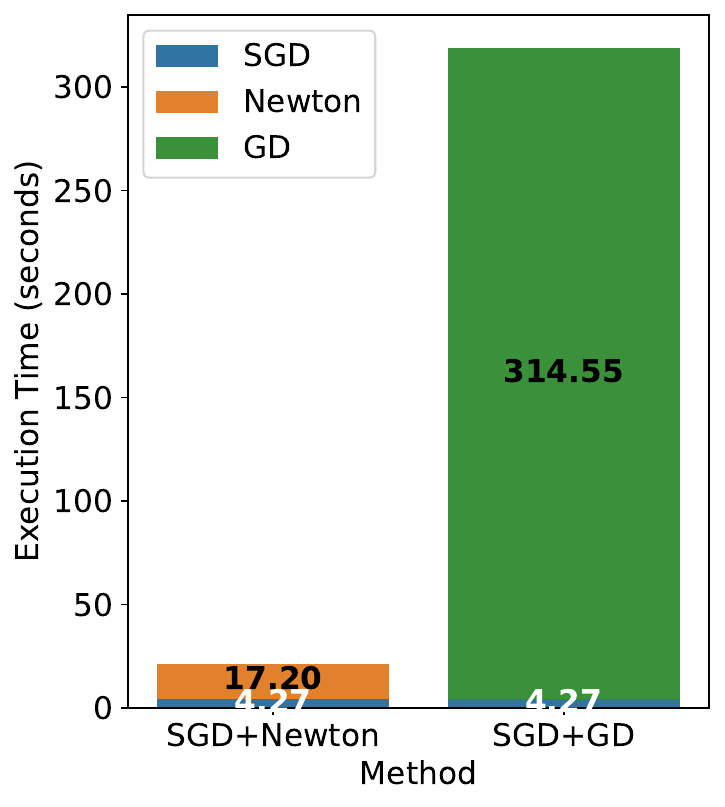}

\caption{3D Point Cloud Registration. 
}
\label{fig:3d_pointcloud}
\vspace{-0.7cm}
\end{figure}

 As shown in Figure \ref{fig:3d_pointcloud}, both methods converge, but at different speeds. The relaxed Newton's method converges in 9 iterations with a runtime of 17.20 seconds, while the GD-only method takes 922 iterations (runtime 314.55 seconds) to reach a comparable loss. Additionally, the relaxed Newton's method achieves about 0.6 orders of magnitude improvement in accuracy in terms of $L_2$ error.

\section{Conclusion}\label{sec:conclusion}
In this work, we computed first-order and second-order derivatives for the parameterized regularized optimal transport (OT) distance. Specifically, we derived explicit analytical expressions for the gradient of the Sinkhorn distance and the Hessian of the entropy-regularized OT (EOT) distance with respect to the source data $\mathbf{Y}$.

To address the numerical instability and high memory consumption typically associated with Hessian computation in large-scale, multi-dimensional problems, we developed a fast, stable, and memory-efficient algorithm using spectral analysis of the ill-posed linear system. Our algorithm demonstrated significant improvements in both efficiency and accuracy on benchmark datasets.

Future work may explore further accelerating the computation by sparse approximation of the Hessian in the computation of both the OT distance and its second-order derivatives, spectral analysis of the H-matrix with random datasets with a focus on the asymptotic behavior in the large $N$ or small $\epsilon$ limit, and robust second-order differentiation for general regularized and constrained optimal transport problems.

\section*{Data Availability}
The software package implementing the proposed algorithms can be found on https://github.com/yexf308/OTT-Hessian.

\section*{Acknowledgements}
X.~Li is grateful for partial support by the NSF Award DMS-1847770 and the 2023 UNC Charlotte faculty research grant. F.~Lu is grateful for partial support by NSF DMS-2238486. M.~Tao is grateful for partial support by the NSF Award DMS-1847802, the Cullen-Peck Scholar Award and the Emory-GT AI Humanity Award. F.~Ye is grateful for partial support by Simons Foundation Award MPS-TSM-00002666.

\appendix

\section{Proofs for Spectral Analysis}
\subsection{Proof of Lemma \ref{lemma:normalize_Pi}}

\begin{proof}
First, we observe that
 \begin{align*}
     \text{diag}(\mm \nu)^{-1} \mm\Pi^\top \text{diag}(\mm  \mu)^{-1} \mm\Pi \mathbb{1}_N =    \text{diag}(\mm \nu)^{-1} \mm\Pi^\top \text{diag}(\mm \mu)^{-1} \mm \mu = \text{diag}(\mm \nu)^{-1} \mm\Pi^\top \mathbb{1}_M =\mathbb{1}_N
 \end{align*}
 This shows the eigenvector of $\lambda=1$ is $\mathbb{1}_N$. 

 Furthermore, we show that $\lambda=1$ has multiplicity one by showing that any solution $\mm v$ to  
  $ \mm v=\text{diag}(\mm\nu)^{-1} \mm\Pi^\top \text{diag}(\mm\mu)^{-1} \mm\Pi \mm v $, must have identical entries. Suppose that there exists  a $\mm v= (v_1,\ldots, v_N)$ has non-constant entries. Let $\widetilde {\mm v}$ denote $ \widetilde {\mm v}= \text{diag}(\mm\mu)^{-1} \mm\Pi {\mm v}$. Then, we have $ \max_{1\leq i\leq M}  \widetilde v_i < \max_{1\leq j\leq N} v_j$, because the entries of $\mm\Pi$ are all positive, $\mm\Pi \mathbb{1}_N =\mu$ and each entry of $ \text{diag}(\mm\mu)^{-1} \mm\Pi {\mm v}$ is a weighted average of $v$. Similarly, we have $ \max_{1\leq j\leq N} [ \text{diag}(\mm\nu)^{-1} \mm\Pi^\top  \widetilde{\mm v}] _j < \max_{1\leq i\leq M}  \widetilde v_i $. Since $\mm v=  \text{diag}(\mm\nu)^{-1} \mm\Pi^\top  \widetilde{\mm v}$ , we obtain a contradiction:  $\max_{1\leq j\leq N} v_j < \max_{1\leq j\leq N} v_j$. Hence,  $\mm v= (v_1,\ldots, v_N)$ must have constant entries, the sames as $\mathbb{1}_N $ up to a scalar factor. 
       
        Next, we show that all eigenvalues must be no greater than one. Suppose that $\lambda$ is an eigenvalue with $v$ being its normalized eigenvector. Thus, denoting $u_1=  \mm\Pi\text{diag}(\mm\nu)^{-1}  v$ and $u_2 = \text{diag}(\mm\mu)^{-1} \mm\Pi {\mm v}$, we have
        \[ \lambda  = \frac{{\mm v}^\top  \text{diag}(\mm\nu)^{-1} \mm\Pi^\top \text{diag}(\mm\mu)^{-1} \mm\Pi {\mm v}}{{\mm v}^\top \mm v} =  \innerp{u_1, u_2} \leq \|u_1\|\|u_2\|\leq 1, \]
        where the last inequality follows from the facts that $\|u_1\|^2= \| \mm\Pi\text{diag}(\mm\nu)^{-1}  \mm v\|_2^2 \leq 1$, and that $\|u_2\| =\| \text{diag}(\mm\mu)^{-1} \mm\Pi {\mm v}\|\leq 1$. 
        
                A similar argument shows that $\text{diag}(\mm\mu)^{-1} \mm\Pi \text{diag}(\mm\nu)^{-1} \mm\Pi^\top$ must have simple eigenvalue  $\lambda=1$ and eigenvector $\mathbb{1}_M$.  
\end{proof}

\subsection{Proof of Proposition \ref{prop:condi_Sinkhorn}}

The proof of Proposition \ref{prop:condi_Sinkhorn} is based on an application of Weyl's inequality to study the eigenvalues of the $\mathbf{H}$-matrix under perturbation. The result is of general interest beyond this study, so we state it as a lemma.

\begin{lemma}[Eigenvalues under perturbation]
\label{lemma:eig_perturb}
Let $\mm\Pi,\mm\Pi^*$ be two positive coupling matrices with $\mathbf{A} = \mm\Pi-\mm\Pi^*$ satisfying 
\begin{equation}\label{eq:perturb_bd}
\max_{i}|\sum_{j}A_{ij} |\le \delta_1, \quad \max_{j}|\sum_i A_{ij}|\le \delta_1, \quad \sum_{i,j}  A_{ij}^2\leq \delta_2^2. 
\end{equation}
Then, the eigenvalues of their $\mathbf{H}$ matrices are close:
\[
|\lambda_k(\mathbf{H})- \lambda_k(\mathbf{H}^*) |\leq  \delta_1+\delta_2, \quad 1\leq k\leq N+M. 
\] 
\end{lemma}
\begin{proof}[\textbf{Proof of Lemma \ref{lemma:eig_perturb}}] Note that we can write the 
\begin{align}
    \mathbf{H}- \mathbf{H}^* = \bcm \mathrm{diag}(\mathbf{A}\mathbb{1})  & \mathbf{A} \\ 
     \mathbf{A}^\top &    \mathrm{diag}(\mathbf{A}^\top \mathbb{1})  \ecm =:  \mathbf{E}. 
\end{align}
By Weyl's inequality, we have 
$|\lambda_k(\mathbf{H}^*) - \lambda_k(\mathbf{H})|\le \|\mathbf{E}\|_{op}$. Thus, it suffices to estimate $\|\mathbf{E}\|_{op}$. Note that first that using $|A_{ij}|\leq 1$ and \eqref{eq:perturb_bd}, we have 
$\sum_{i} | \sum_{j} A_{ij} u_i |^2 = \sum_{i} | \sum_{j} A_{ij} |^2 u_i^2 \leq \sum_{i} | \sum_{j} A_{ij} | u_i^2  \leq \delta_1^2 \| \mathbf{u}\|^2, $
 and similarly, $\sum_{j} | \sum_{i} A_{ij} v_j |^2 \leq \delta_1^2 \|\mathbf{v}\|^2 $;  also, 
$\|\mathbf{A}\mathbf{v}\|^2= \sum_{i=1} | \sum_{j}A_{ij}v_j|^2 \leq  \sum_{i=1} [ \sum_{j}A_{ij}^2 \sum_{j}|v_j|^2] \leq \delta_2^2\|\mathbf{v}\|^2,$ and similarly, $\| \mathbf{A}^\top \mathbf{u} \|^2\leq \delta_2^2 \|\mathbf{u}\|^2$. 

Using these four bounds, we have $\|\mathbf{E}\|_{op}^2 = \sup_{\mathbf{u}\in \R^M, \mathbf{v}\in \R^N, \|\mathbf{u}\|^2+\|\mathbf{v}\|^2 =1} \|\mathbf{E} \bcm  \mathbf{ u} \\  \mathbf{ v} \ecm \|^2 $, then
$\|\mathbf{E}\|_{op}^2 = \sum_{i} | \sum_{j} A_{ij} u_i |^2+ \|\mathbf{Av}\|^2 + \sum_{j} | \sum_{i} A_{ij} v_j |^2  + \|\mathbf{A^\top u}\|^2 \leq  \delta_1^2+ \delta_2^2 \leq \delta_1+\delta_2$.  
Combining with Weyl's inequality, we conclude the proof. 
\end{proof}

\begin{proof}[\textbf{Proof of Proposition \ref{prop:condi_Sinkhorn}}]
The bound for the eigenvalues in \eqref{eq:eig_perturb} follows from \eqref{eq:assum_perturb} and Lemma \ref{lemma:eig_perturb}. 
To prove the bounds for the condition number, we have shown that $\lambda_1: = \lambda_1(\mathbf{H}^*)= \frac{1}{N}+\frac{1}{M}$ and 
$ \lambda_{N+M-1}:=\lambda_{N+M-1}(\mathbf{H}^*)  = \frac{1}{2}\lambda_1\left[1- \sqrt{1- 4\Delta}\right] \in [\lambda_1\Delta, 2\lambda_1\Delta]
$ 
 with $\Delta = \lambda_1^{-2} \rho(\mm\Pi^\top \mm\Pi) $.  For $\delta+\delta_2 =t \lambda_1 \Delta$ with $t\in[0,1)$, Eq.\eqref{eq:eig_perturb} implies that   
\begin{align*}
(1-t)\lambda_1\Delta\leq   \lambda_{N+M-1} -\delta-\delta_2
 &\leq    \lambda_{N+M-1}(\mathbf{H})\leq \lambda_{N+M-1} + \delta+\delta_2 \\
 &\leq (2+t)\lambda_1\Delta \\
\lambda_1(1-t\Delta) \leq   \lambda_{1} -(\delta+\delta_2)
& \leq    \lambda_{1}(\mathbf{H})\leq \lambda_1 + \delta+\delta_2 \leq \lambda_1(1+t\Delta).  
\end{align*}
Consequently, we obtain the bounds by noting that  
\begin{align*}
\frac{1+t\Delta}{(2+t)\Delta} \leq \frac{\lambda_1 - (\delta+\delta_2)}{\lambda_{N+M-1} + (\delta+\delta_2) } 
\leq \kappa(\mathbf{H}) 
\leq   \frac{\lambda_1 + (\delta+\delta_2)}{\lambda_{N+M-1} - (\delta+\delta_2) }  \leq \frac{1+t\Delta}{(1-t)\Delta}.  
\end{align*}
\end{proof}

\section{Equally spaced points on the unit circle}

In this section, we will derive the precise asymptotic behavior of the spectral properties of $\mathbf{H}^*$ as a function of the regularization strength $\epsilon$ and the sample size $N$ in Example 4.8.

\begin{proof} 
 
\textbf{Part (a)}:  Note that the cost matrix is 
	$C_{ij} = \|{y}_i - {y}_{j}\|_2^2 
	= 4\sin^2\left(x_{|j-i|}/2\right)$
and  the Gibbs kernel $K_{ij}=\exp \left(-\frac{4\sin^2\left(|i-j|\pi/N\right)}{\epsilon} \right)$ is circulant, whose rows and columns sum is $\lambda_1(\mathbf{K})=\sum_{j=0}^{N-1}\exp\left(-\frac{4\sin^2\left(j\pi /N\right)}{\epsilon}\right)$. Then the matrix $\mm\Pi^*=\frac{\mathbf{K}}{\lambda_1(\mathbf{K}) N} $ satisfies the uniform marginal constraints on $\mm\Pi^*$, hence, it is the optimal coupling matrix due to uniqueness of the solution of the constraint optimization.

\textbf{Part (b)} follows directly from Corollary 4.6.  

 \textbf{Part (c)}, We first compute the largest two eigenvalues of $\mathbf{K} $. Recall that the matrix $\mathbf{K}$ is symmetric and positive-definite, so its singular values are the same as its eigenvalues.
Since $\mathbf{K}$ is circulant, the first two eigenvalues of $\mathbf{K}$ are (see e.g., \cite{gray2006toeplitz}), 
we have 
\begin{align*}
\lambda_1(\mathbf{K})  = \sum_{j=0}^{N-1} \exp\left(-\frac{4\sin^2(\frac{j\pi}{N})}{\epsilon}\right), \lambda_2(\mathbf{K})=\sum_{j=0}^{N-1} \exp\left(-\frac{4\sin^2(\frac{j\pi}{N})}{\epsilon}\right)\cos\left(\frac{2j\pi}{N} \right) 
\end{align*}
Meanwhile, combining Part (a) and Part (b), we have $\lambda_{2N-1}(\mathbf{H}^*)= \frac{1}{N} (1 - \frac{\lambda_2(\mathbf{K})}{\lambda_1(\mathbf{K})})$. Thus, to study the limits, we first study the limit of $\lambda_1(\mathbf{K})$ and $\lambda_2(\mathbf{K})$. 
	
As $N\rightarrow +\infty$,  the Riemann summations in $\lambda_1$ and $\lambda_2$ approaches the integrals  
\begin{align*}  \lim_{N\rightarrow +\infty} \frac{\lambda_1(\mathbf{K})}{N} 
          &= \frac{1}{\pi}\int_{0}^{\pi} \exp\left(-\frac{4\sin^2(x)}{\epsilon}\right)dx 
           =\exp\left(-\frac{2}{\epsilon}\right) I_0\left(\frac{2}{\epsilon}\right) \\  \lim_{N\rightarrow +\infty} \frac{\lambda_1(\mathbf{K})}{N} 
          &= \frac{1}{\pi}\int_{0}^{\pi} \exp\left(-\frac{4\sin^2(x)}{\epsilon}\right)dx 
           =\exp\left(-\frac{2}{\epsilon}\right) I_0\left(\frac{2}{\epsilon}\right) \\
     \lim_{N\rightarrow +\infty} \frac{\lambda_2(\mathbf{K})}{N}  
           &= \frac{1}{\pi}\int_{0}^\pi \exp\left(-\frac{4\sin^2(x)}{\epsilon}\right)\cos(2x )dx 
           =\exp\left(-\frac{2}{\epsilon}\right) I_1\left(\frac{2}{\epsilon}\right), 
\end{align*}
       where $I_1(x)$ and $I_2(x)$ are the modified Bessel functions of first kind. Then, the limit of the second smallest eigenvalue of $\mathbf{H}^*$ is 
 \begin{align*}
       \lim_{N\rightarrow+\infty}(N\cdot\lambda_{2N-1}(\mathbf{H}^*))&=1-\lim_{N\rightarrow+\infty}\frac{\lambda_2(\mathbf{K})}{\lambda_1(\mathbf{K})} = 1-\frac{I_1(2/\epsilon)}{I_0(2/\epsilon)}.
\end{align*}
  When $\epsilon $ is small, we can expand $I_1(x)$ and $I_0(x)$ around $x=+\infty$,  
          \begin{align*}
        \lim_{N\rightarrow+\infty}(N\cdot\lambda_{2N-1}(\mathbf{H}^*)) =  1-\frac{\frac{(\epsilon/2)^{1/2}}{\sqrt{2\pi}}-\frac{3(\epsilon/2)^{3/2}}{8\sqrt{2\pi}}+O(\epsilon/2)^{5/2}}{\frac{(\epsilon/2)^{1/2}}{\sqrt{2\pi}}+\frac{(\epsilon/2)^{3/2}}{8\sqrt{2\pi}}+O(\epsilon/2)^{5/2}} = \frac{\epsilon}{4} + O(\epsilon^2). 
          \end{align*}
        Then, we have 
        $\lim_{\epsilon\rightarrow 0^+}\lim_{N\rightarrow +\infty} \epsilon \frac{\lambda_1(\mathbf{H}^*) }{\lambda_{2N-1}(\mathbf{H}^*) }=8 $,  which gives \eqref{eq:limit_N_then_e}.  

To prove \eqref{eq:limit_e_then_N}, note first that when $N$ is fixed and $\epsilon \rightarrow 0^+$, $\lambda_1(\mathbf{K})$ and $\lambda_2(\mathbf{K})$ are approximated by the three largest terms, 
 \begin{align*}
    \lambda_1(\mathbf{K})  =  1 + 2\exp(-\frac{4}{\epsilon}\sin^2(\frac{\pi}{N})) + O_\epsilon, \lambda_2(\mathbf{K}) 
       = 1 + 2\exp(-\frac{4}{\epsilon}\sin^2(\frac{\pi}{N}))\cos(\frac{2\pi}{N}) + O_\epsilon
  \end{align*}
  where $O_\epsilon := O\left( \exp\left(-\frac{4\sin^2(2\pi/N)}{\epsilon}\right) \right)$. 
  Consequently,   
     \begin{align*}
        & \exp\left(\frac{4}{\epsilon }\sin^2(\frac{\pi}{N})\right)N\lambda_{2N-1}(\mathbf{H}^*) =   \frac{2 (1- \cos(\frac{2\pi}{N})) +  O_\epsilon}{1 + 2\exp\left(-\frac{4}{\epsilon}\sin^2(\frac{\pi}{N})\right)+ O_\epsilon}. 
     \end{align*}
Taking the limits with $\lim_{N\rightarrow+\infty}\lim_{\epsilon\rightarrow 0^+}$, we obtain  
 \begin{align*} \lim_{N\rightarrow+\infty}\lim_{\epsilon\rightarrow 0^+}N^3\exp\left(\frac{4\sin^2(\pi/N)}{\epsilon }\right)\lambda_{2N-1}(\mathbf{H}^*) = 4\pi^2.
  \end{align*}
\end{proof}

\bibliographystyle{plain}
\bibliography{ref_unsupervised,latentDS,ref_doubly_sto_Matrix,ref_regularization23_09}

\end{document}